\documentclass[times]{nmeauth}
\usepackage[numbers]{natbib}
\usepackage{amsmath}
\usepackage{amscd}
\usepackage{amssymb}
\usepackage{graphicx}
\usepackage{url}
\usepackage{color}
\usepackage{listings}
\usepackage{booktabs}
\usepackage{caption}
\usepackage{textcomp}
\usepackage{mathtools}
\usepackage{mcode}
\usepackage{hyperref}
\usepackage{hypcap}
\usepackage{float,afterpage}
\usepackage{subfigure}
\usepackage[ruled,vlined]{algorithm2e}
\definecolor{DarkBlue}{rgb}{0.00,0.00,0.55}
\definecolor{Black}{rgb}{0.00,0.00,0.00}
\hypersetup{
    linkcolor = DarkBlue,
    anchorcolor = DarkBlue,
    citecolor = DarkBlue,
    filecolor = DarkBlue,
    urlcolor = DarkBlue,
    colorlinks  = true,
}
\graphicspath{{./figures-final/}}
\newcommand{\dx}{\ \mathrm{d}x}
\def\volumeyear{2015}

\begin{document}
\runningheads{P.~E.~Farrell and C.~Maurini}{Numerical solution of a variational fracture model}
\title{Linear and nonlinear solvers for variational phase-field models of brittle fracture}

\author{P.~E.~Farrell\affil{1,2} and C.~Maurini\affil{3},\corrauth}
\address{
\affilnum{1} Mathematical Institute, University of Oxford, Oxford, UK \break
\affilnum{2} Center for Biomedical Computing, Simula Research Laboratory, Oslo, Norway \break
\affilnum{3} Institut Jean Le Rond d'Alembert,  Sorbonne Universit\'es, UPMC, Univ Paris 06, CNRS, UMR 7190, France
}

\corraddr{Corrado Maurini, Institut Jean Le Rond d'Alembert,  Sorbonne Universit\'es, UPMC, Univ Paris 06, CNRS, UMR 7190, France.  E-mail: corrado.maurini@upmc.fr}

\cgs{A Center of Excellence grant from the Research Council of Norway to the Center for Biomedical Computing at Simula Research Laboratory.}
\cgsn{Engineering and Physical Sciences Research Council (UK); Agence Nationale de la Recherche (France)}{EP/K030930/1 and ANR-13-JS09-0009}

\begin{abstract}
The variational approach to fracture is effective for simulating the nucleation
and propagation of complex crack patterns, but is computationally demanding.
The model is a strongly nonlinear non-convex variational inequality that demands
the resolution of small length scales.  The current standard algorithm for its
solution, alternate minimization, is robust but converges slowly and demands the
solution of large, ill-conditioned linear subproblems.  In this paper, we
propose several advances in the numerical solution of this model that
improve its computational efficiency.  We reformulate alternate
minimization as a nonlinear Gauss-Seidel iteration and employ over-relaxation to
accelerate its convergence; we compose this accelerated alternate minimization
with Newton's method, to further reduce the time to solution; and we formulate
efficient preconditioners for the solution of the linear subproblems arising in
both alternate minimization and in Newton's method. We investigate the
improvements in efficiency on several examples from the literature; the new
solver is 5--6$\times$ faster on a majority of the test cases considered.
\end{abstract}

\keywords{fracture, damage, variational methods,  phase-field, nonlinear Gauss-Seidel, Newton's method}

\maketitle


\section{Introduction} \label{sec:introduction}

Cracks may be regarded as surfaces where the displacement field may be
discontinuous.  Fracture mechanics studies the nucleation and propagation of
cracks inside a solid structure.  Variational formulations recast this
fundamental and difficult problem of solid mechanics as an optimization problem.
The variational framework naturally leads to regularized phase-field
formulations based on a smeared description of the discontinuities.  These
methods are attracting an increasing interest in computational mechanics. The
aim of our work is to propose several improvements in the linear and nonlinear
solvers used in this framework.

The code for the algorithms proposed in this paper, and the thermal shock
example of section \ref{sec:tshock}, are included as supplementary material\footnote{The code is available online at \url{https://bitbucket.org/pefarrell/varfrac-solvers}.}.
 
\subsection{Variational formulation of fracture and gradient damage models}
The variational approach to fracture proposed by Francfort and Marigo
\cite{FraMar98} formulates  brittle fracture as the minimization of an energy
functional that is the sum of the elastic energy of the cracked solid and the energy
dissipated in the crack.  The simplest fracture mechanics model, due to Griffith
\cite{Gri21}, assumes that the cracked solid $\Omega\setminus\Gamma$ is linear
elastic and that the surface energy is proportional to the measure of the
cracked surface $\Gamma$.  The crack energy per unit area is the fracture
toughness $G_c$, a material constant. In this case the energy functional to be
minimized is
\begin{equation}
E(u,\Gamma) = \int_{\Omega\setminus\Gamma} \dfrac{1}{2}A_0 \varepsilon(u)\cdot\varepsilon(u)\dx+ G_c \mathcal{S}(\Gamma),
\end{equation}
where $u$ is a vector-valued displacement field,
$\varepsilon(u)=\mathrm{sym}(\nabla u)$ is the second order tensor associated to
the linearised strains inside the material, $A_0$ the fourth order elasticity
tensor of the uncracked solid, and $\mathcal{S}(\Gamma)$ is the Hausdorff
surface measure of the crack set $\Gamma$. In a quasi-static time-discrete
setting, given an initial crack set $\Gamma_0$, the cracked stated of the solid
can be found by incrementally solving the following unilateral minimization
problem \cite{FraBouMar08,FraMar98}:
\begin{equation}
\mathrm{arg\,min}\{E(u,\Gamma), u\in C_{\bar u}(\Omega\setminus\Gamma), \Gamma\supset \Gamma_{i-1}\},
\label{FDP}
\end{equation}
where 
\begin{equation}
\mathcal{C}_{\bar u}(\Omega)\equiv \{ u\in H^1(\Omega, \mathcal R^n), u=\bar u\text{ on }\partial_{\bar u}\Omega\}
\end{equation}
is the space of admissible displacements, $\partial_{\bar u}\Omega$ is the part
of the boundary where the Dirichlet conditions are prescribed and
$H^1(\Omega\setminus\Gamma,  R^n)$ denotes the Sobolev space of vector fields
defined on $\Omega\setminus\Gamma$ with values in $R^n$.  The minimization
problem above is labelled unilateral because the  crack set cannot decrease in
time.  This problem is quasi-static and rate-independent, so that time enters
only via the irreversibility constraint.  The numerical solution of the
\emph{free-discontinuity} problem \cite{AmbFusPal00} above is prohibitive,
because of the difficulty related to the discretization of the unknown crack set
$\Gamma$ where the displacement may jump.

To bypass this issue, Bourdin et al.  \cite{BouFraMar00} transposed to fracture
mechanics a regularization strategy introduced by Ambrosio and Tortorelli
\cite{AmbTor92} for free-discontinuity problems arising in image segmentation
\cite{MumSha89}.  The regularized model approaches the solution of \eqref{FDP}
by the solution of 
\begin{equation} \label{eqn:mainproblem}
\mathrm{arg\,min}\{\mathcal E(u,\alpha), u\in \mathcal{C}_{\bar u}(\Omega), \alpha\in \mathcal{D}_{\alpha_{i-1}}\},
\end{equation}
with the regularized energy functional 
 \begin{equation}
 \mathcal{E}(u, \alpha)=\int_\Omega \left[ {\dfrac{1}{2} a(\alpha)A_0\varepsilon(u)\cdot\varepsilon(u) } +\dfrac{G_c}{c_w}\left (\frac{w(\alpha)}{\ell} +  \ell\, \nabla\alpha\cdot \nabla\alpha\right)\right]\dx, \label{regenergy}
\end{equation}
with $c_w=4 \int_0^1\sqrt{w(\alpha)}\,\mathrm{d} \alpha$.  In this formulation
$\alpha$ is a smooth scalar field, that can be interpreted as  damage, and
$\ell$ is an additional parameter controlling the localization of $\alpha$.
With $\alpha_{i-1}$ the solution at the previous time step and denoting by $\partial_{\bar \alpha}\Omega$  the part
of the boundary where the Dirichlet conditions are prescribed on $\alpha$, the admissible space
for $\alpha$ is a convex cone imposing the unilateral box constraint
\begin{equation}
\mathcal{D}_{\bar\alpha}(\Omega)\equiv \{ u\in H^1(\Omega, \mathcal R), \,
\bar{\alpha} \leq \alpha \leq 1\text{ a.e. in }\Omega,\,\alpha=\bar{\alpha} \text{ on }\partial_{\bar \alpha}\Omega\},
\end{equation}
which prevents self-healing. Following \cite{AmbTor92}, Bourdin et
al.~\cite{BouFraMar00} uses $a(\alpha)=(1-\alpha)^2+k_\ell$ and
$w(\alpha)=\alpha^2$, with $k_{\ell}=o(\ell)$.  With these conditions it is
possible to show through asymptotic methods ($\Gamma$-convergence) that the
solutions of the \emph{global} minimization problem \eqref{eqn:mainproblem} tend to the   solutions of the \emph{global} minimization problem \eqref{FDP} as $\ell \to 0$
\cite{Gia05}.  In the regularized problem, the $\alpha$-field localizes in bands
of thickness on the order of $\ell$ giving smeared representation of the cracks
which is energetically equivalent to the Griffith model (the dissipated energy
per unit crack surface is $G_c$).  This  behaviour is preserved for a large
class of functions $w$ and $a$.  

Similar ``smeared'' crack models have been developed in other contexts.  In the
physics community, they are regarded as phase-field approximations developed by
adapting the Ginzburg-Landau theory of phase transitions \cite{HakKar09}.  In
mechanics, they are regarded as gradient damage models
\cite{LorGod11,PhaAmoMar11}. The $\alpha$ field is a damage variable that
modulates the elastic stiffness and introduces an energy dissipation. In this
context, \eqref{regenergy} can be regarded as a model \emph{per se} for the
evolution of damage in the material, and one can associate a physical meaning to
the internal length and to evolutions following local minima of the functional
\eqref{regenergy}.  In particular, one can show that in quasi-static evolutions
ruled by a \emph{local minimality} condition, the internal length is regarded as
a constitutive parameter that controls the critical stress in the material
before failure.  We refer the reader to \cite{PhaAmoMar11} and Section
\ref{sec:examples} for further details on this point. 

Whilst the original model of Bourdin et al.~\cite{BouFraMar00} assumes small
deformations, isotropic materials,  quasi-static evolution, and allows for
interpenetration of crack lips in compression, recent contributions include
extensions to dynamics \cite{BouLarRic11,BorVerSco12,SchWilKuh14}, multiphysics
couplings \cite{AbdAri12,WheWicWol14}, anisotropic materials \cite{LiPecMil14},
large elastic deformations \cite{DelLanMar07,ClaKna14,HesWei14},  cohesive fracture \cite{VerBor13},  and compressive
failure with unilateral contact at crack lips
\cite{AmoMarMau09,LanRoy09,FreRoy10,MieWelHof10,AmbGerDe-15}, plates and shells \cite{AmiMilShe14,AmbLor16}, thin films \cite{LeoBabBou14}. Other works
\cite{PhaAmoMar11,LorGod11,PhaMarMau11,PhaMar13b} discuss how the choice of the
functions $a$ and $w$ in \eqref{regenergy} affects the properties of the
solutions, by analytical and numerical investigations. On the basis of these
results, recent numerical works \cite[see e.g.][]{BouMarMau14} adopt the choice
$w(\alpha)=\alpha$ and $a(\alpha)=(1-\alpha)^2 + k_{\ell}$, which we also employ
in the rest of this paper.  We refer the reader to \cite{PhaAmoMar11} for a
comparative analysis of this model and the original model in Bourdin et
al.~\cite{BouFraMar00}.

In the remainder of this paper, we discuss the numerical solution of the minimization problem
\eqref{eqn:mainproblem}, after a standard finite element discretisation. We focus on the simplest model, neglecting the effect of geometrical non-linearities and the non-symmetric behaviour of fracture in traction and compression. More complex physical effects drastically modify the character of the numerical problems to be solved and require further problem-specific developments that are outside the scope of this work. 

\subsection{The optimization problem and current algorithms}
\label{sec:currentalg}
The minimization problem \eqref{eqn:mainproblem} problem is numerically challenging, for the following reasons:
\begin{enumerate}
\item the functional is non-convex and thus the minimization problem in general
admits many local minimizers;

\item the irreversibility of damage, required to have a thermodynamically
consistent model and to forbid crack self-healing, introduces bound constraints
on the damage variable $\alpha$ and demands the solution of variational inequalities;

\item the problem size after discretization is usually very large, because the minimizers of
\eqref{regenergy} are typically characterized by localization of damage and
elastic deformations in bands of width on the order of $\ell$. This width is
usually very small with respect to the simulation domain, and the mesh size
must be small enough to resolve the bands;

\item the linear systems to be solved are usually very badly conditioned,
because of the presence of damage localizations where the elastic stiffness
varies rapidly from the undamaged value to zero.
\end{enumerate}

At each loading step, the minimization of \eqref{eqn:mainproblem} is an optimization
problem with necessary optimality conditions: find
$(u,\alpha)\in\mathcal{C}_{\bar u}\times\mathcal{D}_{\bar \alpha_{i-1}}$
satisfying the first order optimality conditions
\begin{eqnarray} \label{eqn:kkt}
\mathcal{E}_u(u, \alpha; v) = 0, \forall v \in \mathcal{C}_0,\qquad  \mathcal{E}_{\alpha}(u, \alpha; \beta) \geq 0, \forall \beta \in \mathcal{D}_{0}
\end{eqnarray}
with 
\begin{align}
\mathcal{E}_u(u, \alpha; v) &= \int_\Omega  {a(\alpha)A_0\varepsilon(u) : \varepsilon(v) }\dx,\\
\mathcal{E}_{\alpha}(u, \alpha; \beta) &= \int_\Omega\left( \dfrac{a'(\alpha)}{2} A_0 \varepsilon(u) : \varepsilon(u) \beta + \dfrac{G_c}{c_w}\frac{w'(\alpha)}{\ell}\beta+ 2 \dfrac{G_c}{c_w} \ell \nabla\alpha\cdot \nabla\beta\right)\dx,
\end{align}
where $\mathcal{E}_u$ and $\mathcal{E}_{\alpha}$ are the Fr\'echet derivatives
of the energy with respect to $u$ and $\alpha$.  Because of the unilateral
constraint on the damage field $\alpha$, the first order optimality conditions
on $\alpha$ form a variational inequality.
The linearization of these conditions is: 
find $(\hat{u}, \hat{\alpha}) \in \mathcal{C}_0 \times \mathcal{D}_0$ such that for all $(v, \beta) \in \mathcal{C}_0 \times \mathcal{D}_0$
\begin{align} \label{eqn:newton}
\mathcal{E}_{uu}(u, \alpha; v, \hat{u}) + \mathcal{E}_{u \alpha}(u, \alpha; v, \hat{\alpha}) &= -\mathcal{E}_u(u, \alpha; v), \nonumber \\
\mathcal{E}_{\alpha u}(u, \alpha; \beta, \hat{u}) + \mathcal{E}_{\alpha
\alpha}(u, \alpha; \beta, \hat{\alpha}) &\geq -\mathcal{E}_\alpha(u, \alpha; \beta),
\end{align}
where
\begin{subequations}
\begin{align}
\mathcal{E}_{uu}(u, \alpha, v, \hat{u}) &= \int_\Omega {a(\alpha)A_0\varepsilon(v)
: \varepsilon(\hat u) }\dx, \label{eqn:elasticity} \\
\mathcal{E}_{u \alpha}(u, \alpha; v, \hat{\alpha}) &= \mathcal{E}_{\alpha u}(u,
\alpha; \hat{\alpha}, v) = \int_\Omega {a'(\alpha)A_0\varepsilon(u) : \varepsilon(v)
\hat{\alpha}}\dx, \label{eqn:offdiag} \\
\mathcal{E}_{\alpha \alpha}(u, \alpha; \beta, \hat{\alpha}) &= \int_\Omega\left[
\left(\dfrac{a''(\alpha)}{2} A_0\varepsilon(u) : \varepsilon(u)  +\dfrac{G_c}{c_w}\frac{w''(\alpha)}{\ell} \right) \beta \hat{\alpha} + 2 \dfrac{G_c}{c_w} \ell \nabla
\beta \cdot \nabla\hat{\alpha}\right]\dx. \label{eqn:helmholtz}
\end{align}
\end{subequations}
Note that the bilinear form $\mathcal{E}_{uu}$ is akin to a standard linear
elasticity problem, and that the bilinear form $\mathcal{E}_{\alpha \alpha}$
is akin to a Helmholtz problem.

The most popular algorithm for the solution of the system \eqref{eqn:kkt} is the
alternate minimization method proposed by Bourdin et al.~\cite{BouFraMar00}.
This algorithm rests on the observation that while the minimization problem
\eqref{eqn:mainproblem} is nonconvex, the functional is convex separately in $u$ or
$\alpha$ if the other variable is fixed. Alternating minimization consists of
alternately fixing $u$ and $\alpha$ and solving the resulting smaller
minimization problem, iterating until convergence. At each iteration before
convergence the optimization subproblem has a unique solution with a lower
energy, and thus the algorithm converges monotonically to a stationary point
\cite{Bou07a}.  The algorithm is detailed in Algorithm \ref{alg:am}.

\begin{algorithm}[htp]
\SetAlgoLined
\KwResult{A stationary point of \eqref{eqn:mainproblem}.}
 Given $(u_{i-1}, \alpha_{i-1})$, the state at the previous loading step. \\
 Set $(u^0, \alpha^0) = (u_{i-1}, \alpha_{i-1})$.

 \While{not converged}{
  Find $u^k\in\mathcal{C}_{\bar u}: \mathcal{E}_u(u, \alpha^{k-1}; v) = 0 \quad \forall v\in \mathcal{C}_{0}$. \\
  Find $\alpha^k \in \mathcal{D}_{\alpha_{i-1}}: \mathcal{E}_{\alpha}(u^k, \alpha; \beta) \geq 0,\quad \forall \beta\in \mathcal{D}_{0}$.
 }

 Set $(u_i, \alpha_i) = (u^k, \alpha^k)$.
 \caption{Standard alternate minimization}
 \label{alg:am}
\end{algorithm}

The first subproblem, finding the updated displacement given a fixed damage,
involves solving a standard linear elasticity problem, but with a strongly
spatially varying stiffness parameter. The second subproblem, finding the
updated damage given a fixed displacement, involves solving a variational
inequality where the Jacobian is a generalized Helmholtz problem, again with
spatially varying coefficients. The standard termination criterion used in
\cite{BouFraMar00} is to stop when the change in the damage field drops below a
certain tolerance. Another approach \cite{AmbGerDe-15} is to stop based on a
normalized change in the energy. Miehe et al.~\cite{MieHofWel10a}  perform a single alternate minimization iteration and propose the use of an adaptive time-stepping.

The main drawback of alternate minimization is its slow convergence rate.  This
motivates the development of alternative approaches using variants of Newton's
method for variational inequalities \cite{LorGod11}, such as active set
\cite{facchinei1998} or semismooth Newton methods \cite{ulbrich2011}. Newton's
method is quadratically convergent close to a solution, but its convergence is
erratic when a poor initial guess is supplied  \cite{LorBad04,HeiWheWic15}.
Numerical experience indicates that Newton's method alone does not converge
unless extremely small continuation steps are taken. Recent attempts to address these convergence issues include the use of
continuation methods \cite{SinVerBor16} or globalization devices such line searches and trust regions \cite{GerDe-}. Moreover,
Newton-type method result in a large system of linear equations to be solved at
each iteration; in prior work direct methods have been employed, limiting the
scalability of the approach.

In this work we make several contributions to the solution of
\eqref{eqn:mainproblem}. First, we cheaply accelerate alternate
minimization by interpreting it as a nonlinear Gauss-Seidel method and applying
over-relaxation. Second, we further reduce the time to solution by composing
alternate minimization with an active set Newton's method, in such a way that
inherits the robustness of alternate minimization and the asymptotic quadratic
convergence of Newton-type methods.  Third, we design scalable linear solvers
for both the alternate minimization subproblems, and the coupled Jacobians of
the form \eqref{eqn:newton} arising in active-set Newton methods.

Following Bourdin et al.\ \cite{BouFraMar00}, we spatially discretize the
problem with standard piecewise linear finite elements on unstructured
simplicial meshes \cite{Neg99}.  This discretization converges to local
minimizers of the Ambrosio--Tortorelli functional \cite{BelCos94}.  Alternatives
proposed in the literature, but not considered here, includes isogeometric
approaches \cite{BorVerSco12}. Adaptive remeshing is a valuable method of
improving to computational efficiency \cite{BurOrtSul10}. The presence of thin
localisation bands renders anisotropic remeshing strategies \cite{ArtForMic15}
particularly attractive. Our work on the linear and non-linear solver is
potentially synergetic with these other efforts to improve computational
efficiency. 

The paper is organised as follows.  Section \ref{sec:nonlinear_solvers} presents
our improved nonlinear solver. The underlying linear solvers and preconditioners
are discussed in section \ref{sec:linear_solvers}. Section \ref{sec:examples}
introduces three fundamental test problems that we use to assess the performance
of the solvers. The results of the corresponding numerical experiments are
reported in section~\ref{sec:numerical_experiments}. Finally,
we conclude in section~\ref{sec:conclusions}.

\section{Nonlinear solvers} \label{sec:nonlinear_solvers}

In this section we propose several improvements to the nonlinear solver employed
for the minimization of the regularized energy functional \eqref{eqn:mainproblem}. The
first improvement is to reinterpret alternate minimization as a nonlinear
Gauss-Seidel iteration: this naturally suggests employing an \emph{over-relaxed}
Gauss-Seidel approach, which we discuss in section \ref{sec:oram}. This
over-relaxation greatly reduces the number of iterations required for
convergence, with minimal computational overhead. The second improvement is to
use alternate minimization as a preconditioner for Newton's method
\cite{brune2013}, as discussed in section \ref{sec:oramn}. By 
combining these, our solver enjoys the robust convergence of (over-relaxed)
alternate minimization and the rapid convergence of Newton's method. Alternate
minimization is used to drive the approximation within the basin of convergence
of Newton's method; once this is achieved, Newton's method takes over and solves
the nonlinear problem to convergence in a handful of iterations. As our
numerical experiments in section \ref{sec:numerical_experiments} demonstrate, this strategy
is faster than relying on alternate minimization alone, even with
over-relaxation.

\subsection{Over-relaxed alternate minimization} \label{sec:oram}
\begin{algorithm}[htp]
\SetAlgoLined
\KwResult{A stationary point of \eqref{eqn:mainproblem}.}
 Given $(u_{i-1}, \alpha_{i-1})$, the state at the previous loading step, and
 the over-relaxation parameter $\omega \in (0, 2)$. \\
 Set $(u^0, \alpha^0) = (u_{i-1}, \alpha_{i-1})$.

 \While{not converged}{
  Find $\tilde{u}^k\in\mathcal{C}_{\bar u}: \mathcal{E}_u(u, \alpha^{k-1}; v) = 0 \quad \forall v\in \mathcal{C}_{0}$. \\
  Set $\delta u^k = \tilde{u}^k - u^{k-1}$. \\
  Set $u^{k} = u^{k-1} + \omega \delta u^k$.
  \vspace{0.3cm}

  Find $\tilde{\alpha}^k\in\mathcal{D}_{\alpha_{i-1}}: \mathcal{E}_{\alpha}(u^{k}, \tilde{\alpha}^k; \beta) \geq 0 \quad \forall \beta\in \mathcal{D}_{0}$. \\
  Set $\delta \alpha^k = \tilde{\alpha}^k - \alpha^{k-1}$. \\
  Choose the largest $\bar{\omega} \in (0, \omega)$ so that $\alpha^{k-1} +
  \bar{\omega} \delta \alpha^k \in \mathcal{D}_{\alpha_{i-1}}$. \\
  Set $\alpha^{k} = \alpha^{k-1} + \bar{\omega} \delta \alpha^k$.
 }

 Set $(u_i, \alpha_i) = (u^k, \alpha^k)$.
 \caption{Over-relaxed alternate minimization (ORAM)}
 \label{alg:oram}
\end{algorithm}

In the block-Gauss-Seidel relaxation method for linear systems, the solution variables
are partitioned; at each iteration, some variables are frozen and a linear
subproblem is solved for the remaining free variables; the updated values for
these variables are used in the solution of the next subset.  Similarly, a
nonlinear block-Gauss-Seidel relaxation first solves a nonlinear subproblem for one
subset of the variables, then uses those updated values to solve for the next
subset, and so on \cite{ortega1970}. Alternate
minimization is precisely a nonlinear block-Gauss-Seidel method that iterates
between the displacement and damage variables. Just as over-relaxation can
accelerate linear Gauss-Seidel \cite{young1950}, it can also
accelerate nonlinear Gauss-Seidel \cite{ortega1970}. Therefore, we
augment the standard alternate minimization algorithm with a simple
over-relaxation approach, Algorithm \ref{alg:oram}. 
The state before and after each alternate minimization substep are compared to
determine the update direction, and over-relaxation is applied along that
direction with relaxation parameter $\omega$. In the damage step, the bound
constraint on $\alpha$ is enforced during the line search: if a step with
$\omega$ would be infeasible, the algorithm sets $\bar{\omega}$ to the midpoint of
$[1, \omega]$, and repeats this recursively until the update to $\alpha$ retains
feasibility. (The question of infeasibility does not arise for $\omega < 1$.)

The literature on over-relaxation methods is vast, and we briefly summarise the
main points here.  In linear successive over-relaxation (SOR) applied to a
matrix $A$, the convergence depends on the spectral radius of the SOR iteration
matrix
\begin{equation*}
M = (D - \omega L)^{-1}(\omega U + (1 - \omega)D),
\end{equation*}
where $D, -L$ and $-U$ are the diagonal, lower triangular and upper triangular
components of $A$. Essentially, over-relaxation attempts to choose an $\omega$
that reduces $\rho(M)$. A similar result holds for block SOR \cite{arms1956}.
Kahan \cite{kahan1958} proved that $\omega \in (0, 2)$ is a necessary condition
for the convergence of SOR, i.e.~for $\rho(M) < 1$. Ostrowski
\cite{ostrowski1954} proved that this is sufficient for convergence in the case
where $A$ is symmetric and positive-definite.  For nonlinear SOR, Ortega and
Rheinboldt \cite[Theorem 10.3.5]{ortega1970} proved the surprising result that
the asymptotic convergence rate depends on the spectral radius of the SOR
iteration matrix evaluated at the Jacobian of the residual evaluated at the
solution. Nonlinear Gauss-Seidel methods ($\omega = 1$) can also be extended to
minimisation problems with constraints, under the name of block coordinate
descent.  We are not aware of any analysis of over-relaxation in the context of
constraints, or in the infinite dimensional setting, as would be necessary here
for a proof of convergence; however, the numerical experiments of section
\ref{sec:numerical_experiments} demonstrate that convergence was achieved for
all problems with all values of $\omega \in (0, 2)$ attempted, and that
over-relaxation can significantly reduce the number of iterations required for
convergence on difficult problems.

\subsection{Choosing the relaxation parameter $\omega$}
The number of iterations required depends sensitively on the choice of
$\omega$. Extrapolating from the nonlinear SOR theory, we hypothesize
that the optimal $\omega$ is that that minimizes the spectral radius
of the SOR iteration matrix associated with the unconstrained degrees
of freedom at the minimizer. Unfortunately, identifying this $\omega$
\emph{a priori} appears to be difficult: such an analysis would rely
on the spectral properties of the Hessian at the unknown minimizer
\cite{reid1966}, which are not in general known. In this work we rely
on the na\"ive strategy of numerical experimentation on coarser problems,
and defer an automated scheme for choosing $\omega$ to future work.

\subsection{Composing over-relaxed alternate minimization with Newton} \label{sec:oramn}
Even with over-relaxation, achieving tight convergence of the optimization
problem takes an impractical number of iterations (on the order of hundreds or
thousands for difficult problems).  Therefore, instead of driving the
optimization problem to convergence with ORAM, we use it instead to bring the
iteration within the basin of convergence of a Newton-type method, Algorithm
\ref{alg:oramn}. There are two main problems to solve in designing such a
composite solver: first, deciding when to switch from ORAM to Newton, and
second, handling the possible failure of the Newton-type method.

The inner termination criterion for the over-relaxed alternate minimization used
in this work is based on the norm of the residual of the optimality conditions
\eqref{eqn:kkt}. As the optimality conditions are a variational inequality, it
is not sufficient to merely evaluate a norm of $[\mathcal{E}_u,
\mathcal{E}_\alpha]$, because the feasibility condition should enter in the
termination criterion. Instead, the residual is defined via a so-called
nonlinear complementarity problem (NCP) function $\Phi$: a function which is
zero if and only if the variational inequality is satisfied \cite{munson2000}.
In this work we use the Fischer--Burmeister NCP-function, which is described in
section \ref{sec:vi}. A typical inner termination criterion might be to
switch when the norm of the Fischer--Burmeister residual has decreased by two
orders of magnitude, although the choice taken should vary with the difficulty
of the problem considered. If the inner tolerance is too tight, an excessive
number of alternate minimization iterations will be performed before switching
to Newton; if the inner tolerance is too loose, then the Newton iteration may
not converge and the extra cost of solving Jacobians yields no advantage.

If the Newton-type method diverges (possibly significantly increasing the
residual), it can be handled in one of two ways. The first is to check at the
end of an outer iteration whether Newton's method reduced the residual: if not,
discard the result of Newton and continue with more alternate minimization
iterations. The second is to choose a Newton-type method that is guaranteed to
monotonically decrease the norm of the residual, or to terminate with failure:
this property is achieved by complementing the Newton iteration with a
backtracking line search. This latter option was implemented in our experiments.
If the Newton method fails to achieve a sufficient reduction, the outer
composite solver simply reverts to alternate minimization to bring the solution
closer to the basin of convergence. In this way, the robustness and monotonic
convergence of alternate minimization is combined with the quadratic asymptotic
convergence of Newton's method.

\begin{algorithm}[htp]
\SetAlgoLined
\KwResult{A stationary point of \eqref{eqn:mainproblem}.}
 Given $(u_{i-1}, \alpha_{i-1})$, the state at the previous loading step. \\
 Set $(u^0, \alpha^0) = (u_{i-1}, \alpha_{i-1})$.

 \While{not converged}{
  Set $\|\Phi^0\|$ to be the norm of the residual of the optimality conditions evaluated at $(u^k, \alpha^k)$. \\
  \vspace{0.5cm}

  \While{$\|\Phi^k\| / \|\Phi^0\| > $ inner tolerance}{
    Apply over-relaxed alternate minimization, Algorithm \ref{alg:oram}.
  }
  \vspace{0.5cm}

  \While{not converged and maximum iterations not reached}{
    Apply an active-set Newton method with backtracking line search, such as Algorithm \ref{alg:fs}.
  }
 }

 Set $(u_i, \alpha_i) = (u^k, \alpha^k)$.
 \caption{Over-relaxed alternate minimization combined with Newton (ORAM-N)}
 \label{alg:oramn}
\end{algorithm}

\subsection{Reduced-space active set method} \label{sec:vi}
Both the damage subproblem and the subproblem to be solved at each
coupled Newton iteration are variational inequalities, which when
discretized yield complementarity problems. In this section we briefly
review the Newton-type method used to solve these complementarity problems,
a reduced-space active set method implemented in PETSc \cite{benson2006}.
While semismooth Newton methods have gained significant popularity in recent
years, the reduced-space method employed in this work makes devising
preconditioners for the linear system to be solved in section
\ref{sec:solvingrs} more straightforward.

A mixed complementarity problem (MCP) is defined by a residual $F: \mathbb{R}^n
\rightarrow \mathbb{R}^n$, a lower bound vector $l \in \mathbb{R}_{-\infty}^n$,
and an upper bound vector $x \in \mathbb{R}_{\infty}^n$, where
$\mathbb{R}_{\infty} = \mathbb{R} \cup \{\infty\}$ and $\mathbb{R}_{-\infty} =
\mathbb{R} \cup \{-\infty\}$.  A solution $x \in \{x \in \mathbb{R}^n: l \leq x \leq
u\}$ satisfies $\mathrm{MCP}(F, l, u)$ iff, for each component $i$, precisely
one of the following conditions holds:
\begin{align}
F_i(x) = 0 &\text{ and } l_i < x_i < u_i \nonumber \\
F_i(x) \ge 0 &\text{ and } l_i   =  x_i                    \\
F_i(x) \le 0 &\text{ and } \phantom{l_i =  }\, x_i   =  u_i \nonumber.
\end{align}

A special case of a mixed complementarity problem is the choice $l = 0, u =
\infty$, which is referred to as a nonlinear complementarity problem (NCP),
$\mathrm{NCP}(F)$. For clarity, the algorithm will be described in the context
of NCPs; the extension to MCPs is straightforward \cite{munson2000}.

An NCP-function $\phi: \mathbb{R} \times \mathbb{R} \rightarrow \mathbb{R}$ is a function
with the property that $\phi(a, b)=0 \iff a \ge 0, b \ge 0, ab = 0$, i.e. that $a$
solves $\mathrm{NCP}(b)$. An example is the Fischer--Burmeister function
\cite{fischer1992}
\begin{equation}
\phi_{\mathrm{FB}}(a, b) = \sqrt{a^2 + b^2} - a - b.
\end{equation}
NCP-functions are useful because it is possible to reformulate an NCP as a
rootfinding problem. Given an NCP-function, it is possible to define a residual of $\mathrm{NCP}(F)$:
\begin{equation}
\Phi_i(x) = \phi(x_i, F_i(x)).
\end{equation}
A solution $x$ satisfies $\mathrm{NCP}(F)$ iff $\Phi(x) = 0$.
While $\Phi(x)$ is semismooth, its squared-norm $\|\Phi(x)\|^2$ is smooth
\cite{facchinei1997}.

\begin{algorithm}[htp]
\SetAlgoLined
\KwResult{A solution of $\mathrm{NCP}(F)$.}
 Given $x^0$, the initial guess.

 \While{$\|\Phi(x^k)\| > $ tolerance}{
 Compute the active and inactive sets $A$ and $N$ via \eqref{eqn:reducedis}. \\
 Set $d_A = 0$. \\
 Solve the reduced Newton step \eqref{eqn:reducednewton} for $d_N$. \\
 Choose the step length $\mu$ such that $\|\Phi\|^2$ is minimized, via line search on $\pi\left[x^k + \mu d\right]$;
 if this search direction fails, use the steepest descent direction instead.
 }
 \caption{Reduced-space active set method.}
 \label{alg:fs}
\end{algorithm}

At each iteration, the algorithm constructs a search direction $d$.
The search direction is defined differently for the \emph{active} and \emph{inactive}
components of the state.
Given an iterate $x$ and a fixed zero tolerance $\zeta > 0$, define the active set
\begin{equation} \label{eqn:reducedis}
A = \{i : x_i \le \zeta \text{ and } F_i(x) > 0 \},
\end{equation}
and define the inactive set $N$ as its complement in $\{1, \dots, n\}$. The active set represents
a hypothesis regarding which variables will be zero at the solution.
At each iteration, the active subvector of the search direction is zeroed.
For the inactive component of the search direction, a Newton step is performed.
The inactive component is defined by approximately solving
\begin{equation} \label{eqn:reducednewton}
J_{N, N} d_N = F_N,
\end{equation}
where $J$ is the Jacobian of the residual $F$.
The submatrix retains any symmetry and positive-definiteness properties of the underlying
Jacobian \cite{benson2006}.
Given this search direction, a line search is performed with merit function
$\|\Phi(x)\|^2$, with each candidate projected on to the bounds with projection
operator $\pi$. If this line search fails, the steepest descent direction
is used instead.
The algorithm is listed in Algorithm \ref{alg:fs}. A major advantage of this approach over
other algorithms is that the linear systems to be solved in \eqref{eqn:reducednewton} are of
familiar type: they are submatrices of PDE Jacobians, which have been well studied in the literature.
This familiarity is exploited to design suitable preconditioners for \eqref{eqn:reducednewton} in
section \ref{sec:solvingrs}.

\section{Linear solvers and preconditioners}
\label{sec:linear_solvers}
With this configuration of nonlinear solvers, there are three linear subproblems
to be solved: linear elasticity for the displacement field with fixed damage; a
Helmholtz-like operator for the damage field with fixed displacement; and
(submatrices of) the coupled Jacobian of the optimality conditions. As it is
desirable to solve finely-discretized problems on supercomputers, it is
important to choose scalable iterative solvers and preconditioners for each
subproblem. These are discussed in turn.

\subsection{Linear elastic subproblem}
The linear elastic problem is symmetric and positive-definite, and hence, the
method of conjugate gradients \cite{hestenes1952} is used. However, the problem
is poorly conditioned due to the strong variation in stiffness induced by damage
localization, and appropriate preconditioners must be employed.  The Krylov
solver is preconditioned by the \texttt{GAMG} smoothed aggregation algebraic
multigrid preconditioner \cite{adams2004}, which is known to be extremely
efficient for large-scale elasticity problems.

For algebraic multigrid to be efficient, it is essential to supply the algorithm
with the near-nullspace of the operator, eigenvectors associated with
eigenvalues of small magnitude \cite{falgout2006}. The elasticity problem
without damage has a near-nullspace consisting of the rigid body modes of the
structure; with damage, the localized variation in stiffness induces additional
near-nullspace vectors. Calculations of the smallest eigenmodes of the
elasticity operator with SLEPc \cite{hernandez2005} indicate that if the
structure is partitioned into two or more undamaged regions separated by damaged
regions, the elasticity operator has additional near-nullspace vectors
associated to independent rigid body motions of the separate regions. For
example, suppose algebraic multigrid alone (no Krylov accelerator) is used to
solve the elasticity problem arising with the converged damage field of the
problem of the traction of a bar (section \ref{sec:traction_bar}). With no
nullspace configured, convergence is achieved in 2004 multigrid V-cycles; if
only the entire rigid body modes are supplied, convergence is achieved in 50
V-cycles; and if the additional near-nullspace vectors corresponding to the
partition are supplied, then convergence is achieved in 6 V-cycles.

While these additional near-nullspace vectors assist the convergence of the algebraic multigrid
algorithm, they are very difficult to compute, as they depend on the damage field itself. Therefore
in this work we do not supply these additional near-nullspace vectors, supplying only the rigid body
modes of the entire structure. When a Krylov method is used to accelerate the convergence of the
algebraic multigrid, the ratio of iteration counts between the full and partial near null-spaces
decreases from approximately 10 to approximately 2. However, it may be possible to improve the
convergence of the linear elasticity problem by approximating the additional near-nullspace vectors
arising due to damage. This could be of significant benefit, as this phase constitutes a large
proportion of the solver time.

\subsection{Damage subproblem}
The inactive submatrix of the Helmholtz problem for damage is also solved with
conjugate gradients and the \texttt{ML} smoothed aggregation multigrid algorithm
\cite{vanek1996,gee2006}, with the near-nullspace specified as the constant
vector.

\subsection{The Newton step} \label{sec:solvingrs}
Let the inactive submatrix of the coupled Jacobian be partitioned as
\begin{equation} \label{eqn:inactjac}
J = \begin{bmatrix}
A & B \\
B^T & C
\end{bmatrix},
\end{equation}
where $A$ is the assembly of linear elasticity operator \eqref{eqn:elasticity},
and $B$ and $C$ are the inactive submatrices of the coupling term
\eqref{eqn:offdiag} and the linearised damage operator
\eqref{eqn:helmholtz}.
The fast iterative solution of block matrices has been a major topic of research
in recent years \cite{murphy2000}, with most preconditioners relying on the approximation of a
Schur complement of the operator. It is straightforward to verify that if
$A$ is invertible, then the inverse of a block matrix like \eqref{eqn:inactjac}
can be written as \cite[equation (3.4)]{benzi2005}
\begin{align}
J^{-1} &=
\begin{bmatrix}
A^{-1} & 0 \\
0 & I
\end{bmatrix}
\begin{bmatrix}
I & -B \\
0 & I
\end{bmatrix}
\begin{bmatrix}
A & 0 \\
0 & S^{-1}
\end{bmatrix}
\begin{bmatrix}
I & 0 \\
-B^T & 0
\end{bmatrix}
\begin{bmatrix}
A^{-1} & 0 \\
0 & I
\end{bmatrix} \\
&=
\begin{bmatrix}
A^{-1} + A^{-1}BS^{-1}B^TA^{-1} & -A^{-1}BS^{-1} \\
-S^{-1}B^TA^{-1} & S^{-1}
\end{bmatrix},
\end{align}
where $S = C - B^T A^{-1} B$ is the (dense) Schur complement matrix of $J$ with
respect to $A$.
In this work, we take the simple approximation $S \approx C$, which yields
the preconditioner
\begin{align} \label{eqn:preconditioner}
P^{-1} &= 
\begin{bmatrix}
A^{-1} & 0 \\
0 & I
\end{bmatrix}
\begin{bmatrix}
I & -B \\
0 & I
\end{bmatrix}
\begin{bmatrix}
A & 0 \\
0 & C^{-1}
\end{bmatrix}
\begin{bmatrix}
I & 0 \\
-B^T & 0
\end{bmatrix}
\begin{bmatrix}
A^{-1} & 0 \\
0 & I
\end{bmatrix} \\
&= \begin{bmatrix}
A^{-1} + A^{-1}BC^{-1}B^TA^{-1} & -A^{-1}BC^{-1} \\
-C^{-1}B^TA^{-1} & C^{-1}
\end{bmatrix}
\end{align}
which requires one application of $C^{-1}$ and two applications of $A^{-1}$ per
preconditioner application. This is implemented in PETSc using the symmetric
multiplicative variant of the \texttt{PCFIELDSPLIT} preconditioner
\cite{brown2012,balay2014}. Both inverse actions are approximated by two
V-cycles of algebraic multigrid. MINRES \cite{paige1975} is employed as the
outer Krylov solver, as far from minimizers the Hessian may not be positive
definite.

\section{Test cases} \label{sec:examples}
In this section we introduce three test cases that are used to assess the
performance of the proposed solvers. These test cases will then be used to
assess the performance of the solver in section \ref{sec:numerical_experiments}.
The first investigates temporally smooth propagation of a single crack driven by
appropriately controlled Dirichlet boundary conditions. The second consists of
the uniaxial traction of a bar, testing crack nucleation.  The last considers a
thermal shock problem involving the nucleation and propagation of a complex
pattern. All test cases consider isotropic homogeneous materials. In this
context, the relevant material parameters are the Poisson ratio $\nu$, the
Young's modulus $E$, the fracture toughness $G_c$, and the internal length
$\ell$. One can show that the internal length may be estimated by knowledge of
the limit stress $\sigma_c$ through the relation
\begin{equation}
\ell = \frac{3}{8} \frac{G_c E}{\sigma_c^2}.
\end{equation}
Dimensional analysis shows that without loss of generality both $G_c$ and $E$
can be set to 1, with a suitable rescaling of the loading. Hence, in all
experiments we fix $G_c = E = 1$, and in addition we fix $\nu = 0.3$.

\subsection{Surfing: smooth crack propagation}
The main advantage of the variational regularized approach to fracture analyzed
in this paper is its ability to compute the propagation of cracks along complex
paths, including crack bifurcation, merging, and possible jumping in time and
space. However, it is desirable to test the numerical algorithm in a simpler
situation where a single preexisting crack is expected to propagate smoothly
along a straight path with an assigned velocity $v$. To this end, we consider
the surfing experiment proposed by Hossein et al.\ \cite{HosHsuBou14}. This
consists of a rectangular slab $\Omega = [0, L] \times [-H/2, H/2]$ of length
$L$ and height $H$ with the Dirichlet boundary condition 
\begin{equation}
u(x_1,x_2,t) = \bar{u}(x_1 - L_c - v\,t,x_2)\quad \text{ on } \partial \Omega
\end{equation}
imposed on the whole boundary of the domain. $\bar{u}$ is the
asymptotic Mode-I crack displacement of linear elastic fracture mechanics
\begin{equation}
\bar u = \dfrac{K_I}{2\mu}\sqrt{\frac{r}{2\pi}}\left(\frac{3-\nu}{1+\nu}-\cos \theta\right) 
\left(\cos{({\theta}/{2})} e_1 + \sin({\theta}/{2}) e_2\right),
\label{surfingloading}
\end{equation}
where $(r,\theta)$ are the polar coordinates, $(e_1, e_2)$ are the Cartesian
unit vectors, $\mu$ is the shear modulus, and $L_c$ is the length of the
preexisting crack. The intensity of the loading is controlled by the stress
intensity factor $K_I$. From the theory we expect that the crack 
propagates at the constant speed $v$ along the line $x_2=0$ for
$K_I=K_I^c=\sqrt{G_c\,E}$.  In the numerical experiments we set $K_I/K_I^c =
1.0$, $v = 1$, $L=2$, $H=1$ and $L_c = 0.05$.

\begin{figure}[ht]
\begin{center}
\begin{minipage}{0.5\textwidth}
\centering
 \includegraphics[width=\textwidth]{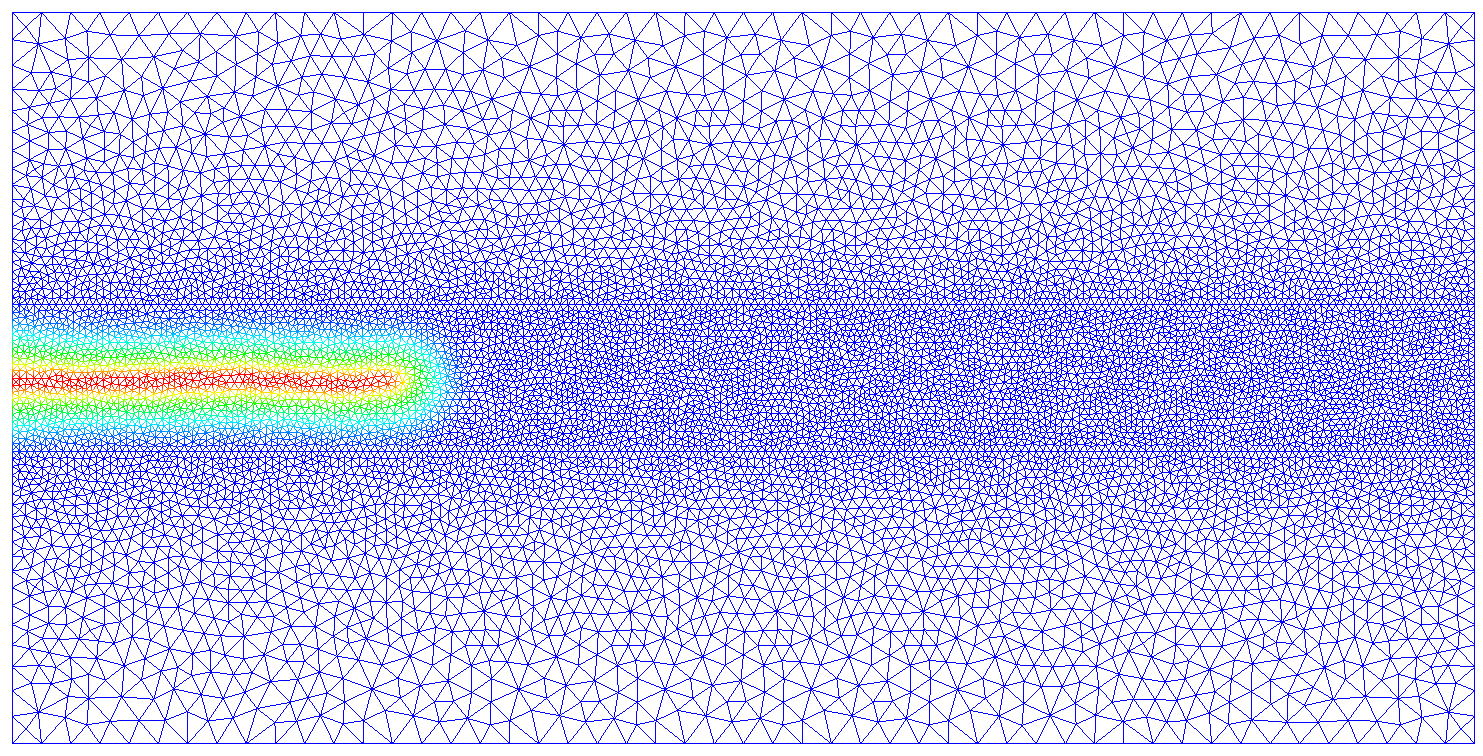} 
\end{minipage}\hfill
\begin{minipage}{0.5\textwidth}
\centering
 \includegraphics[width=\textwidth]{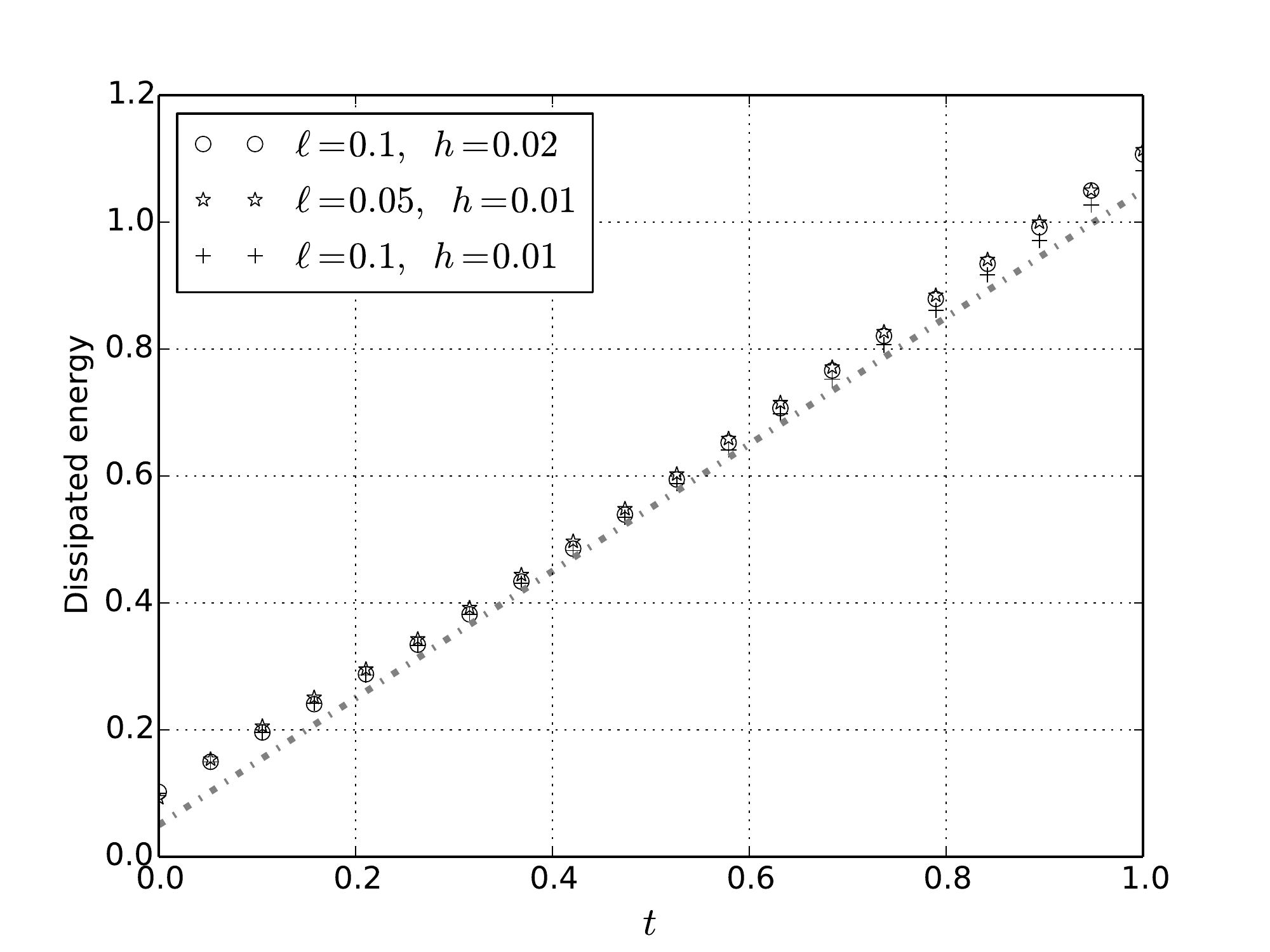} 
\end{minipage}
\end{center}
  \caption{Smooth crack propagation test on a rectangular slab of dimensions
  $2\times 1$ with the surfing loading \eqref{surfingloading} applied on the
  boundary. Left: snapshot of the damage field and mesh for $\ell=0.05$ and
  $h=0.01$. Right: Dissipation energy versus time for  $v=1$ and $K_I=1$
  comparing the results obtained through the damage model when varying the
  internal length  $\ell$ and the mesh size $h$; the continuous line is the
  expected surface energy according to the Griffith model,  corresponding to a
  constant crack speed $v=1$. The given mesh size refers to the typical element
  dimension in the refined band in the middle.}
  \label{fig:surfing}
\end{figure}

Figure~\ref{fig:surfing} reports the results of the corresponding numerical
simulations. This test is particularly useful to verify that the dissipated
energy does not depend on $\ell$ and is equal to the product of the crack length and the
fracture toughness $G_c$. Obviously, in order for this condition to hold,
the discretisation should be changed with the internal length, as $\ell$
controls the width of the localization band. We typically set the mesh size to
$h=\ell/5$.  In the present test, to speed up the benchmarks, we use a
non-uniform mesh respecting this condition only in the band where we expect
the crack to propagate, as shown in Figure~\ref{fig:surfing}. This \emph{a
priori} mesh refinement is exceptional and not applicable in general. In all
other tests, a sufficiently fine uniform mesh will be employed.

\subsection{Traction of a bar} \label{sec:traction_bar}

\begin{figure}[ht]
\begin{center}
  \includegraphics[width=.45\textwidth]{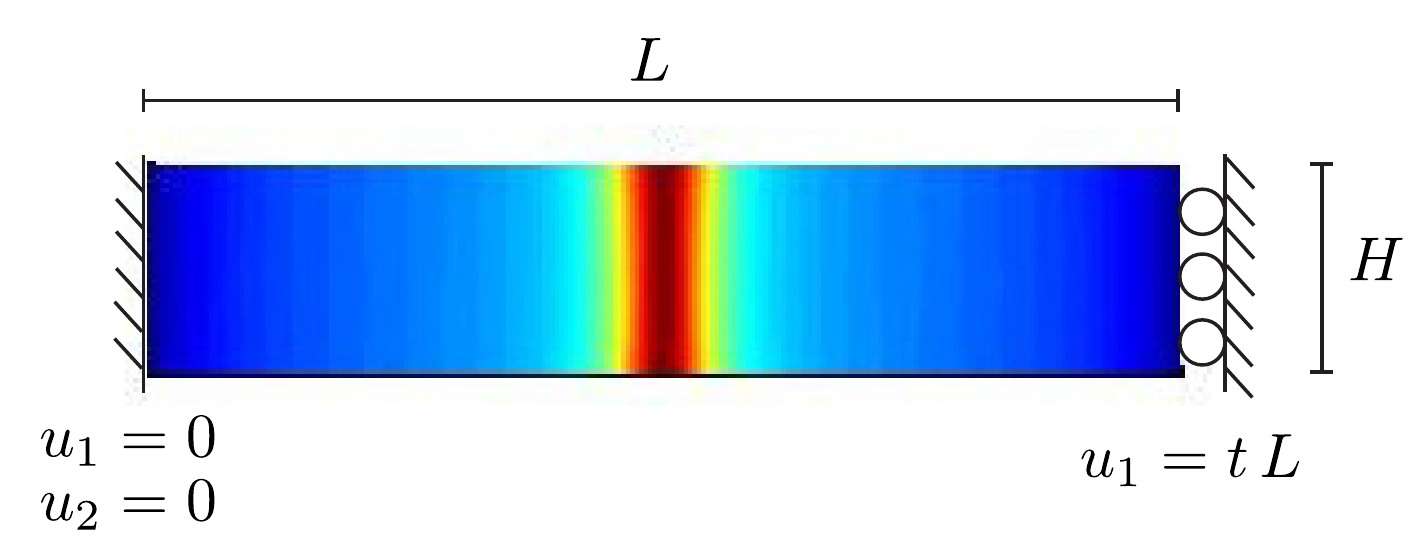} 
  \includegraphics[width=.45\textwidth]{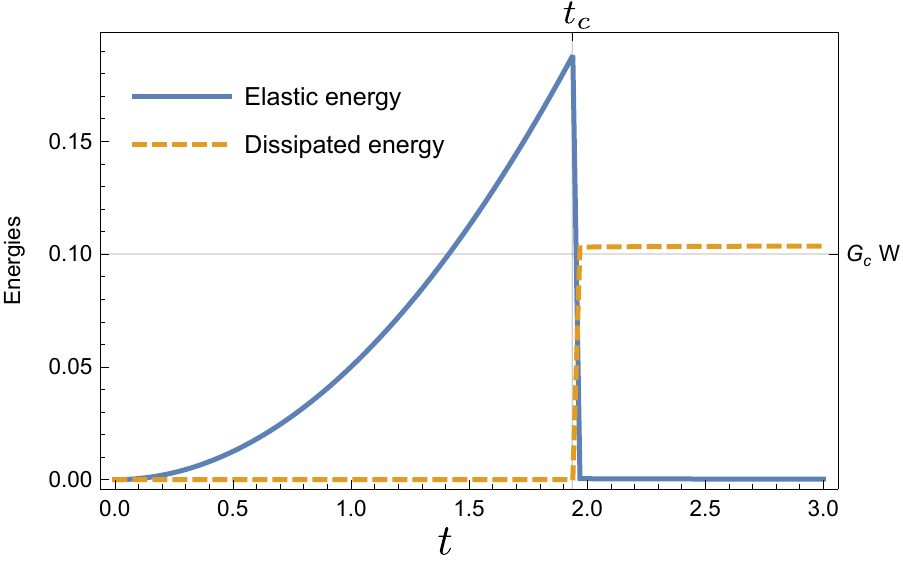} 
\end{center}
  \caption{Uniaxial traction of bar. Left: boundary conditions and damage field
  for $t>t_c$. The width of the localization band may be calculated analytically
  and is given by $2\sqrt{2}\,\ell$. Right: evolution of the energy at the
  solution given by the minimization algorithm as the applied end-displacement
  $t$ is increased.}
  \label{fig:bar}
\end{figure}

A basic problem of fracture mechanics is to estimate the ultimate load before
fracture of a straight bar in traction.  We consider a two dimensional bar of
length $L$  and height $H$ under uniaxial traction with imposed displacement, as
shown in figure \ref{fig:bar}. Analytical studies \cite{PhaMarMau11,PhaAmoMar11}
show that, for $L$ sufficiently greater that $\ell$, a local minimum of the
energy functional \eqref{regenergy} is the purely elastic solution $\alpha=0$
for $t< t_c=\sqrt{{3G_c}/{8 E\ell}}$ and the solution with one crack represented
in figure \ref{fig:bar} for $t>t_c$. The cracked solution has a vanishing
elastic energy and a surface energy given by $G_c\,W$. The test may be easily
extended to a 3D geometry. The critical load $t_c$ is the same in 1D under a
uniaxial stress condition, in 2D plane stress, or in 3D.

\subsection{Thermal shock} \label{sec:tshock}
\begin{figure}[htbp]
\begin{center}
\includegraphics[width=1.\textwidth]{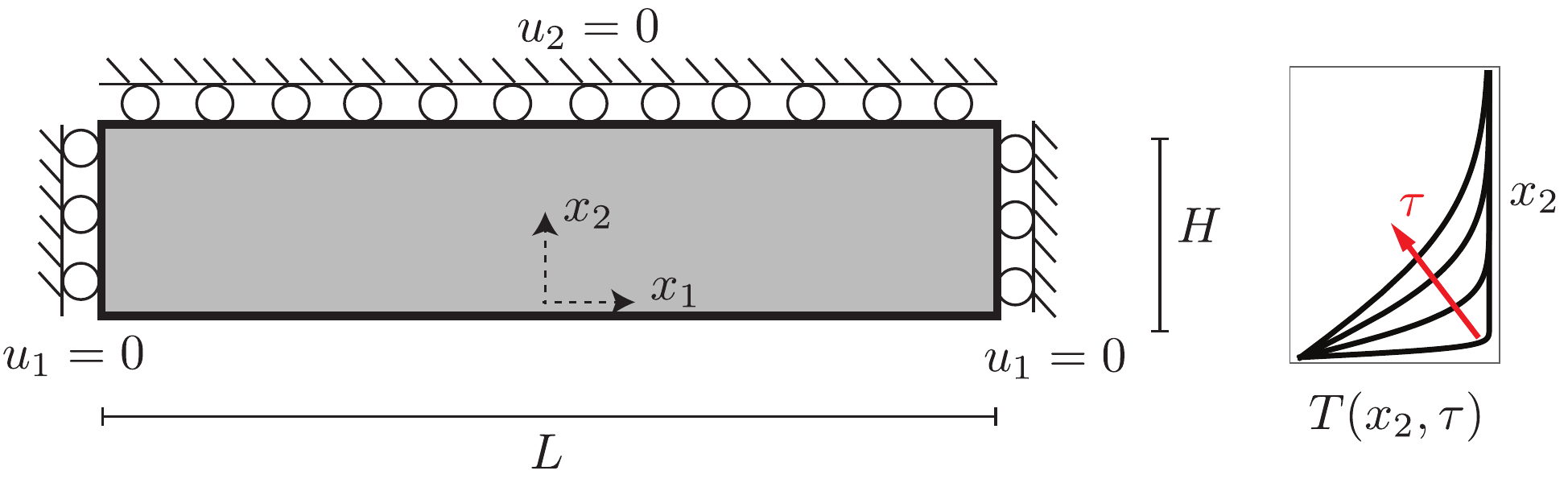}
\caption{Geometry and boundary conditions for the thermal shock problem (left),
where $u_1$ and $u_2$ denotes the two components of the displacement field. The
loading is given by the thermal stress induced by the temperature field
$T(x_2,\tau)$ of \eqref{tempfielderfc}, whose dependence in $x_2$ is sketched on
the right for different times $\tau$. }
\label{tshockgeom}
\end{center}
\end{figure}

The thermal shock problem of a brittle slab  \cite{BouMarMau14} is a challenging
numerical test for the nucleation and propagation of multiple cracks. In
physical experiments \cite{BahWeiMas88,JiaWuLi12}, several ceramic slabs are
bound together, uniformly heated to a high temperature and quenched in a cold
bath, so as to submit the boundary of the domain to a thermal shock. The
inhomogeneous temperature field induces an inhomogeneous stress field inside the
slab, causing the emergence of a complex crack pattern, with an almost periodic
array of cracks nucleating at the boundary and propagating inside the slab with
a period doubling phenomenon.  Following \cite{BouMarMau14,SicMarMau14}, we
consider a simplified model of this experimental test. The computational domain
$\Omega=[-L/2,L/2]\times[0,H]$ (see Figure \ref{tshockgeom}) is a slab of  width
$L$ and height $H$, with a thermal shock applied at the bottom surface $x_2=0$.
At each timestep $\tau_i$, we seek the quasi-static evolution of the cracked
state of the solid  by solving for a stationary point of the following energy
functional:
\begin{equation}
 {	\mathcal{E}_{{\ell}}(u,\alpha)=
		\int_{\Omega} \dfrac{a(\alpha)}{2} A_0 \,\varepsilon_{\mathrm{eff}}(u;\tau):\varepsilon_{\mathrm{eff}}(u;\tau)\dx}		+ \,
		{\dfrac{G_c}{c_w} \int_{\Omega}\left(\dfrac{w(\alpha)}{\ell}  +\,{\ell}\,\nabla \alpha\cdot\nabla\alpha \right)\dx},
\label{energytshock}
\end{equation}
where $\varepsilon_\mathrm{eff}(u;\tau)=\varepsilon(u)-\varepsilon_0(\tau)$ is
the effective elastic deformation accounting for  the thermally induced
inelastic strain
\begin{equation}
\varepsilon_0(\tau)={\beta\, T(\tau)\, I}, \qquad T(\tau) = -\Delta T\,\mathrm{Erfc}\left(\dfrac{x_2/\ell}{\,\tau}\right),
\label{tempfielderfc}
\end{equation}
where $\beta$ is the thermal expansion coefficient.  The temperature field $T$
imposed is the analytical solution of an approximate thermal diffusion problem
with a Dirichlet boundary condition on the temperature for a semi-infinite
homogeneous slab of thermal diffusivity $k_c$. In particular, it neglects the
influence of the cracks on the thermal diffusivity.  The function $
\mathrm{Erfc}$ denotes the complementary error function and $\tau = 2\sqrt{{k}_c
t  /\ell^2}$ a dimensionless time acting as the loading parameter.

As discussed in Bourdin et al.~\cite{BouMarMau14}, the system is governed by
three characteristic lengths: the size of the domain $L$, the internal length
$\ell$, and the Griffith length $\ell_0={G_c}/{E(\beta \Delta T)^2}$. Hence,
choosing $\ell$ as the reference length, the solution depends on two
dimensionless parameters: the mildness of the thermal shock $\ell_0/\ell$ and
the size of the slab $L/\ell$.  Here we perform numerical simulations for fixed
slab dimensions $L = 20$ and internal length $\ell = 1$. We apply the
displacement boundary conditions described in Figure \ref{tshockgeom} and do not
impose any Dirichlet boundary condition on the damage field.  As  can be show by
dimensional analysis, without loss of generality, we set $E=1$, $G_c=1$, $\beta=1$.  We consider the performance of the solver for varying
$\Delta T$ and mesh sizes $h$. 

\begin{figure}[ht!]
  \centering
  \includegraphics[width=\textwidth]{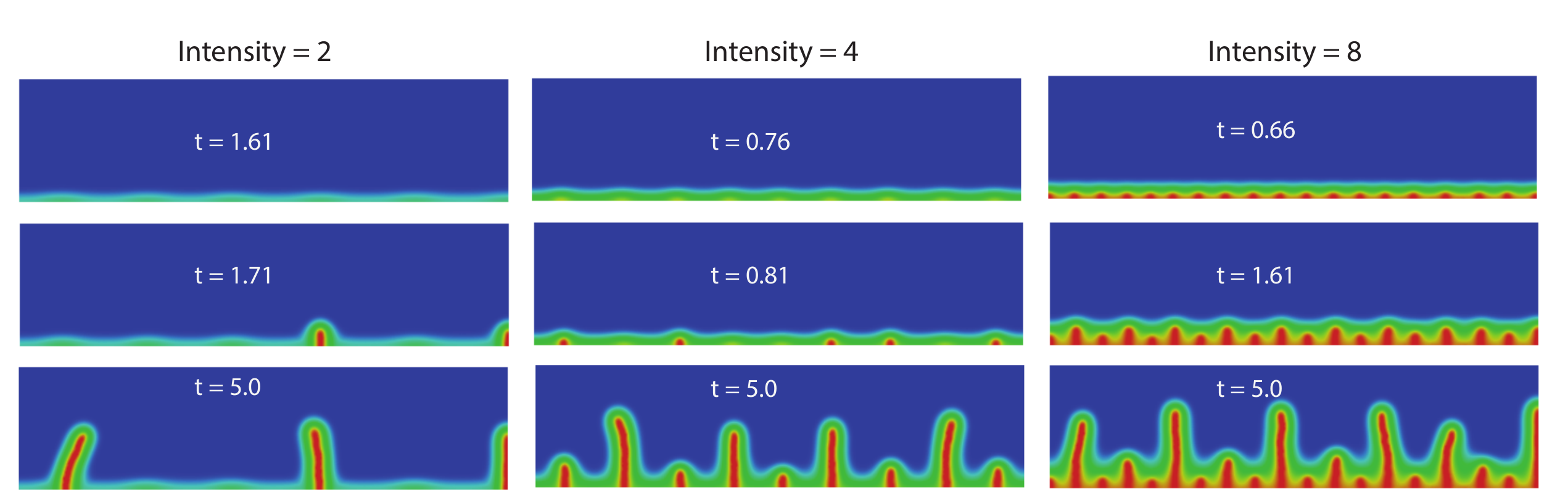}\\
  \caption{Snapshot of the evolution of the damage variable $\alpha$ during the
  evolution (blue: $\alpha=0$; red: $\alpha=1$)  showing the initial solution
  independent of the $x_1$ variable, the emergence of a periodic crack pattern
  and its selective propagation with period doubling. Each column corresponds to
  the result obtained for a specific intensity, increasing from left (2) to
  right (8). Here  $\ell=1$ and the slab  dimensions are $40\times 10$ with a
  mesh size $h=0.2$.}
  \label{fig:tshockevolution}
\end{figure}

\begin{figure}[ht!]
  \centering
  \includegraphics[width=.6\textwidth]{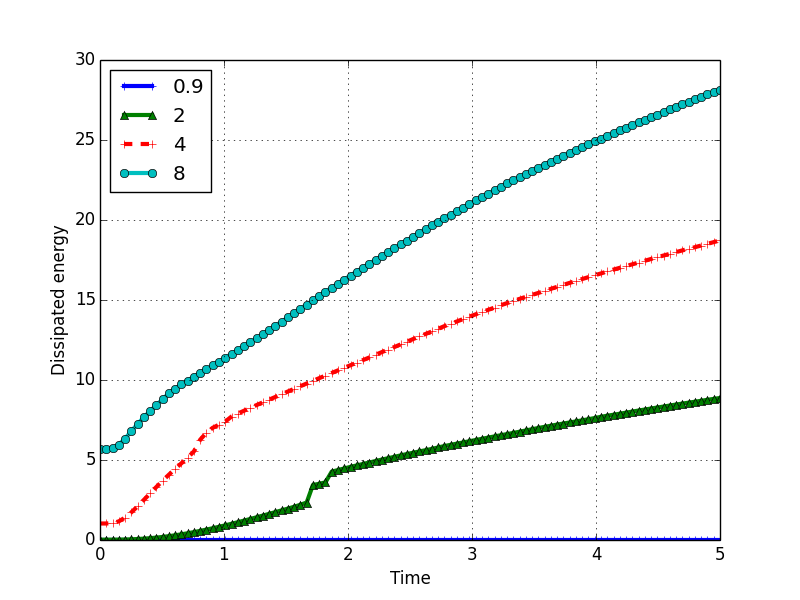}\\
  \caption{Dissipated energy versus time for the thermal shock problem with intensity
  $\Delta T/\Delta T_c$ equal to 0.9, $2$, $4$ and $8$, as in Figure
  \ref{fig:tshockevolution}. For intensity $\Delta T= 2\Delta T_c$, the
  evolution shows two clear jumps in time, corresponding to snapbacks and
  sudden crack growths.  By contrast, the evolution is smooth for  $\Delta T=
  8\Delta T_c$. }
  \label{fig:tshockenergies}
\end{figure}

Analytical and
semi-analytical results are available  for verification purposes and for the
design of the numerical experiments. For $\Delta T<\Delta T_c  =
\sqrt{8 E\ell/3 \beta^2\,G_c}$ the solution is purely elastic with no damage
($\alpha=0$ everywhere) \cite{BouMarMau14,SicMarMau14}.  For $\Delta
T>\Delta T_c $ the solution evolves qualitatively as in Figure
\ref{fig:tshockevolution}, with (i) the immediate creation of an $x$-homogeneous
damage band parallel to the exposed surface, (ii) the bifurcation of this
solution toward an $x$-periodic one, which (iii) further develops in a periodic
array of crack bands orthogonal to the exposed surface. These bands further
propagate with a period doubling phenomenon  (iv).  The three columns in Figure
\ref{fig:tshockevolution} show the phases (ii)-(iv) of the evolution for
$\Delta T/\Delta T_c $ equal to 2, 4 and 8. The wavelength of the oscillations
and the spacing of the cracks increase with  $\Delta T$. In particular
\cite{BouMarMau14} shows that for $\Delta T\gg \Delta T_c$ the initial crack
spacing is proportional to $\sqrt{ \ell_0\ell}$. Figure \ref{fig:tshockenergies}
reports the evolution of the dissipated energy versus time for the three cases
of Figure \ref{fig:tshockevolution}. We note in particular that, while the
evolution is smooth for intense thermal shocks (see the curve $\Delta T = 8
\Delta T_c$), for mild shocks there are jumps in the energy dissipation and
hence in the crack length (see the curve $\Delta T = 2 \Delta T_c$). These jumps
correspond to snap-backs in the evolution problem, where the minimization
algorithm is obliged to search for a new solution, potentially far from the one
at the previous time step.
 
This problem constitutes a relevant and difficult test for the solver.
First, the presence of a large number of cracks renders the
elastic subproblem particularly ill-conditioned, and tests the effectiveness of
the linear subsolvers and the coupled preconditioning strategy. Second, the presence of
bifurcations and snap-backs stresses the convergence and effectiveness of the
outer nonlinear solver. Third, the solution of the overall quasi-static
evolution problem is strongly influenced by the irreversibility condition, testing
the effectiveness of the variational inequality solvers.

\section{Results of numerical experiments}
\label{sec:numerical_experiments}
We present here the results of the numerical experiments that were performed to
assess the performance of the proposed solvers.  All problems were solved to an
absolute $l_2$ residual tolerance of $10^{-7}$. For each test problem, we
analyse the dependence of the results on the relevant physical parameter:
we vary the internal length $\ell$ in the traction and surfing tests,
and the intensity of the loading $\Delta T/\Delta T_c$ in the thermal shock
problem.

\subsection{Over-relaxation} \label{sec:am_acceleration}
We first consider ORAM, the over-relaxation of alternate minimization
described in section \ref{sec:oram}. Each problem of section \ref{sec:examples}
was solved  with values of the
over-relaxation parameter $\omega$ taken from $\{0.8, 1.0, \dots, 
1.8\}$. To consider the effect of over-relaxation alone, the Newton solver
was disabled and all linear solves were performed with LU \cite{amestoy2001}.

\begin{table}
\begin{center}
\begin{tabular}{cccccccc}
\toprule
$\ell$ & $\omega = 0.8$ & $\omega = 1.0$ & $\omega = 1.2$ & $\omega = 1.4$ & $\omega = 1.6$ & $\omega = 1.8$ & reduction \\
\midrule
0.20  & 361 & 233 & 148 & 98 & 174 & 394 & 57.94\% \\
0.10  & 564 & 369 & 251 & 159 & 148 & 307 & 59.89\% \\
0.05  & 1168 & 773 & 537 & 368 & 236 & 320 & 69.47\% \\
0.02  & 2523 & 1680 & 1182 & 835 & 569 & 461 & 72.56\% \\
\bottomrule
\end{tabular}
\caption{Impact of over-relaxation on the surfing case. Standard alternate
minimization converges slowly, and over-relaxation significantly reduces the
number of iterations required.}
\label{tab:am_acceleration_surfing}

\begin{tabular}{cccccccc}
\toprule
$\ell$ & $\omega = 0.8$ & $\omega = 1.0$ & $\omega = 1.2$ & $\omega = 1.4$ & $\omega = 1.6$ & $\omega = 1.8$ & reduction \\
\midrule
0.10  & 124 & 53 & 111 & 181 & 324 & 729 & 0\% \\
0.05  & 120 & 37 & 115 & 185 & 326 & 747 & 0\% \\
0.02  & 132 & 39 & 121 & 195 & 332 & 726 & 0\% \\
0.01  & 139 & 39 & 121 & 186 & 325 & 709 & 0\% \\
\bottomrule
\end{tabular}
\caption{Impact of over-relaxation on the traction case.  Standard alternate
minimization converges rapidly for all values of $\ell$ and over-relaxation
hinders convergence.}
\label{tab:am_acceleration_traction}

\begin{tabular}{cccccccc}
\toprule
$\Delta T/\Delta T_c$ & $\omega = 0.8$ & $\omega = 1.0$ & $\omega = 1.2$ & $\omega = 1.4$ & $\omega = 1.6$ & $\omega = 1.8$ & reduction \\
\midrule
 2  & 3577 & 2364 & 1685 & 1273 & 1095 & 1666 & 53.68\% \\
 4  & 3283 & 2156 & 1548 & 1184 & 1023 & 1611 & 52.55\% \\
 8  & 5100 & 2619 & 1844 & 1354 & 1094 & 1542 & 58.23\% \\
16  & 5097 & 3382 & 2367 & 1669 & 1226 & 1756 & 63.75\% \\
\bottomrule
\end{tabular}
\caption{Impact of over-relaxation on the thermal shock case.
Standard alternate minimization converges slowly, and over-relaxation
significantly reduces the number of iterations required.}
\label{tab:am_acceleration_tshock}
\end{center}
\end{table}

The results for the surfing, traction and thermal shock problems are shown in
Tables \ref{tab:am_acceleration_surfing}, \ref{tab:am_acceleration_traction} and
\ref{tab:am_acceleration_tshock} respectively. In all tables, the reduction
column describes the decrease in iterations for the optimal $\omega$ compared to
standard alternate minimization, $\omega = 1$.  In the traction case, standard
alternate minimization is extremely efficient: a small number of iterations is
required, the number of iterations required does not grow with $\ell$, and
applying any other $\omega$ slows the convergence of the method. By contrast, in
the surfing and thermal shock cases, standard alternate minimisation converges
slowly, and the number of iterations required increases as the physical parameters 
$\ell$ and $\Delta T/\Delta T_c$ are varied. In
this sense, the surfing and thermal shock cases are harder than the traction
case. In these problems, over-relaxation helps significantly, reducing the
number of iterations required by a factor between $1/2$ and $3/4$. Furthermore,
the advantage gained by over-relaxation increases as the problem gets harder.

\subsection{Composition of alternate minimization with Newton's method}
We next consider ORAM-N, the composition of alternate minimization with Newton's method
as described in section \ref{sec:oramn}. For these experiments,
Newton's method was attempted once alternate minimization had reduced the $l_2$
norm of the residual by $10^{-1}$. All linear solves (both for alternate
minimization and Newton's method) were performed with LU, and all alternate
minimizations employed the optimal over-relaxation parameter determined in the
previous experiments ($\omega = 1.6$ for the surfing and thermal shock cases,
$\omega = 1$ for the traction case). The time in seconds was measured for both
approaches, as comparing iteration counts would be irrelevant. The runs
were executed in serial on an otherwise unloaded Intel Xeon E5-4627 3.30 GHz CPU
with 512 GB of RAM.

\begin{table}
\begin{center}
\begin{tabular}{cccc}
\toprule
    & \multicolumn{2}{c}{time (s)} & \\
$\ell$ & alternate minimization alone & composite solver & reduction \\
\midrule
0.20 & 9.11 & 3.68 & 59.60\%\\
0.10 & 27.00 & 15.43 & 42.85\%\\
0.05 & 168.42 & 96.83 & 42.51\%\\
0.02 & 2643.86 & 1886.53 & 28.64\%\\
\bottomrule
\end{tabular}
\caption{Combining alternate minimization with Newton's method for the surfing
case. This further reduces
the runtime of the solver compared to over-relaxed alternate minimization.}
\label{tab:am_newton_surfing}

\begin{tabular}{cccc}
\toprule
    & \multicolumn{2}{c}{time (s)} & \\
$\ell$ & alternate minimization alone & composite solver & reduction \\
\midrule
0.10 & 2.94 & 2.60 & 11.56\%\\
0.05 & 5.94 & 5.71 & 3.87\%\\
0.02 & 41.90 & 29.32 & 30.02\%\\
0.01 & 195.97 & 175.06 & 10.67\%\\
\bottomrule
\end{tabular}
\caption{Combining alternate minimization with Newton's method for the traction
case. In this case, the gains are modest.}
\label{tab:am_newton_traction}

\begin{tabular}{cccc}
\toprule
    & \multicolumn{2}{c}{time (s)} & \\
$\Delta T/\Delta T_c$ & alternate minimization alone & composite solver & reduction \\
\midrule
2 & 312.29 & 215.68 & 30.94\%\\
4 & 319.87 & 220.73 & 30.99\%\\
8 & 335.38 & 238.28 & 28.95\%\\
16 & 386.57 & 287.48 & 25.63\%\\
\bottomrule
\end{tabular}
\caption{Combining alternate minimization with Newton's method for the thermal
shock case. This further
reduces the runtime of the solver compared to over-relaxed alternate
minimization.}
\label{tab:am_newton_tshock}
\end{center}
\end{table}

The results for the surfing, traction and thermal shock problems are shown in
Tables \ref{tab:am_newton_surfing}, \ref{tab:am_newton_traction} and
\ref{tab:am_newton_tshock} respectively. Again, the traction case is unusual
compared to the other two: while the gains are marginal in the traction
case, composition yields a worthwhile and consistent reduction in runtime for the
other tests. If a more robust semismooth Newton solver were available,
the speedup from composition would further increase.

\subsection{Preconditioning the full Jacobian} \label{sec:tshock_preconditioning}
The preconditioner \eqref{eqn:preconditioner} requires inner solvers for the
displacement elasticity operator $A$ and the damage Helmholtz operator $C$. We
first consider the performance of \eqref{eqn:preconditioner} with ideal inner
solvers (LU), to investigate how the iteration counts scale with the physical parameters 
and with mesh size $h$. We then consider the performance with
practical inner solvers, two V-cycles of algebraic multigrid for $A$ and $C$.
Each Jacobian solve was terminated when the $l_2$ norm of the
residual was reduced by a factor of $10^{-6}$, although adaptive tolerance
selection should be used in practical calculations to retain quadratic
convergence of the inexact Newton method \cite{eisenstat1996}. For each
configuration of physical parameters and $h$, the total number of Krylov
iterations required for convergence over all loading steps was divided by the
total number of Newton iterations, to compute the average number of Krylov
iterations required to solve a Jacobian. In these experiments the alternate
minimization was terminated with a relative residual reduction of $10^{-3}$, or
if the absolute residual norm reached $10^{-6}$. As the gains from employing
Newton's method in the traction case were marginal, we consider here only
the surfing and thermal shock problems.

\begin{table}
\begin{center}
\begin{tabular}{cccc}
\toprule
    & \multicolumn{3}{c}{average Krylov iterations} \\
$\ell$ & $h = \ell/5$ & $h = \ell/10$ & $h = \ell/15$ \\
\midrule
0.20 & 6.17 & 9.47 & 12.03\\
0.10 & 8.92 & 10.91 & 13.53\\
0.05 & 10.74 & 13.07 & 15.58\\
0.02 & 13.07 & 14.95 & 18.03\\
\bottomrule
\end{tabular}
\caption{The average Krylov iterations per Newton step for different internal
 lengths $\ell$ and mesh sizes $h$ for the surfing case, with
ideal inner solvers (LU). The preconditioner depends weakly on mesh refinement
and on $\ell$.}
\label{tab:newton_ideal_surfing}

\begin{tabular}{cccc}
\toprule
    & \multicolumn{3}{c}{average Krylov iterations} \\
$\ell$ & $h = \ell/5$ & $h = \ell/10$ & $h = \ell/15$ \\
\midrule
0.20 & 6.33 & 10.57 & 11.64\\
0.10 & 8.50 & 10.93 & 13.68\\
0.05 & 11.13 & 13.53 & 16.15\\
0.02 & 13.53 & 15.61 & 18.80\\
\bottomrule
\end{tabular}
\caption{The average Krylov iterations per Newton step for different internal
 lengths $\ell$ and mesh sizes $h$ for the surfing case, with
practical inner solvers (AMG). Switching to practical inner solvers hardly
affects the convergence of the preconditioner.}
\label{tab:newton_krylov_surfing}

\begin{tabular}{cccc}
\toprule
    & \multicolumn{3}{c}{average Krylov iterations} \\
$\Delta T/\Delta T_c$ & $h = \ell/4$ & $h = \ell/8$ & $h = \ell/16$ \\
\midrule
2 & 10.11 & 7.87 & 6.87\\
4 & 10.56 & 11.44 & 14.70\\
8 & 15.97 & 19.64 & 17.75\\
16 & 23.15 & 22.37 & 22.03\\
\bottomrule
\end{tabular}
\caption{The average Krylov iterations per Newton step for different 
intensities $\Delta T/\Delta T_c$  
and mesh sizes $h$ for the thermal shock case, with
ideal inner solvers (LU). The preconditioner depends weakly on mesh refinement
and on Griffith length.}
\label{tab:newton_ideal_tshock}

\begin{tabular}{cccc}
\toprule
    & \multicolumn{3}{c}{average Krylov iterations} \\
$\Delta T/\Delta T_c$ & $h = \ell/4$ & $h = \ell/8$ & $h = \ell/16$ \\
\midrule
2 & 15.25 & 12.54 & 11.82\\
4 & 12.95 & 13.42 & 18.25\\
8 & 19.08 & 19.65 & 22.51\\
16 & 27.79 & 27.22 & 29.54\\
\bottomrule
\end{tabular}
\caption{The average Krylov iterations per Newton step for different intensities $\Delta T/\Delta T_c$ 
and mesh sizes $h$ for the thermal shock case, with
practical inner solvers (AMG). In this case, using practical inner solvers
somewhat degrades the convergence of the preconditioner.}
\label{tab:newton_krylov_tshock}
\end{center}
\end{table}

The results for the surfing case with ideal and practical inner solvers are
shown in Tables \ref{tab:newton_ideal_surfing} and
\ref{tab:newton_krylov_surfing}, and the corresponding results for the thermal
shock case are shown in Tables \ref{tab:newton_ideal_tshock} and
\ref{tab:newton_krylov_tshock}. In the surfing case, the number of iterations
required grows slowly as the mesh is refined, and grows slowly as $\ell$ is
reduced. However, even for the smallest $\ell$ on the finest mesh, the number
of outer Krylov iterations required is modest, and the results barely differ if
the ideal inner solvers are replaced with practical variants. In the thermal
shock case, the number of iterations required stays approximately constant
as the mesh is refined, and grows slowly with  the intensity  $\Delta T/\Delta T_c$. Here, replacing
the ideal inner solvers with practical variants does have a measurable
cost in iteration count; this could be reduced by tuning the parameters
of the algebraic multigrid algorithm employed, or by employing stronger
inner solvers. These results show that the preconditioner
\eqref{eqn:preconditioner} is a practical and efficient solver for the full
coupled Jacobian, whose performance degrades slowly as the difficulty of
the problem is increased.

\section{Conclusion} \label{sec:conclusions}
In this paper we proposed several improvements to the current standard algorithm
for solving variational fracture models. Over-relaxation is extremely cheap and
simple to implement, but can greatly reduce the number of iterations
required for convergence. Composing over-relaxed alternate minimization with Newton-type methods yields a further
decrease in runtime, although at a more significant development cost. Together, these improvements to
alternate minimization reduce the time to solution by a factor of 5--6$\times$ for the surfing and thermal shock
test cases. Lastly, we proposed and tested preconditioners for the
linear subproblems in alternate minimization and the coupled Jacobian arising in the  Newton iterations when solving the whole problem with a monolithic active set method.
These efforts are complementary to other approaches recently proposed in the literature, such as
adaptive remeshing \cite{ArtForMic15}, adaptive time-stepping, continuation algorithms \cite{SinVerBor16}, or refined line-search techniques \cite{GerDe-} that were not considered in this work.
Our  tests focus only on the simplest settings for variational fracture mechanics assuming small deformations and a simple rate-independent  material behaviour. However,  the developed techniques can be readily adapted to more complex contexts, including  hyperelasticity, viscoelasticity,  and  inertial effects.

These results suggest several directions for future research. It would be highly
desirable to develop a convergence analysis of block over-relaxed nonlinear
Gauss-Seidel for variational inequalities, although we do not anticipate this
will yield constructive insight for the choice of the over-relaxation parameter $\omega$. It may be possible
to design $\ell$-robust preconditioners for the coupled Jacobian (where the
convergence is independent of $\ell$) by choosing appropriate $\ell$-dependent
inner products for the displacement and damage function spaces. If the
appropriate Babu\v{s}ka constants are independent of $\ell$, the convergence
will be also \cite{mardal2011}.

In this work we have considered only the simplest fracture model, assuming small deformations and symmetric behaviour in traction and compression. Our developments on over-relaxation and composition of alternate minimization with Newton could be applied with minor modifications to more complex cases, including for example the tension-compression splitting of the elastic energy to account for the non-interpenetration condition on the crack lips \cite{AmoMarMau09,LanRoy09,FreRoy10,MieWelHof10,AmbGerDe-15}. In this case an efficient solver would require the development of suitable preconditioners for the elastic subproblem and the coupled Jacobian.

\bibliographystyle{wileyj}
\bibliography{literature}

\begin{thebibliography}{10}
\providecommand{\url}[1]{\texttt{#1}}
\providecommand{\urlprefix}{URL }
\expandafter\ifx\csname urlstyle\endcsname\relax
  \providecommand{\doi}[1]{doi:\discretionary{}{}{}#1}\else
  \providecommand{\doi}{doi:\discretionary{}{}{}\begingroup
  \urlstyle{rm}\Url}\fi

\bibitem{FraMar98}
Francfort G, Marigo JJ. Revisiting brittle fracture as an energy minimization
  problem. \emph{Journal of the Mechanics and Physics of Solids}  1998;
  \textbf{46}(8):1319--1342.

\bibitem{Gri21}
Griffith A. The phenomena of rupture and flow in solids. \emph{Philosophical
  Transactions of the Royal Society A}  1921; \textbf{221}:163--198.

\bibitem{FraBouMar08}
Francfort G, Bourdin B, Marigo JJ. The variational approach to fracture.
  \emph{Journal of Elasticity}  2008; \textbf{91}(1-3):5--148.

\bibitem{AmbFusPal00}
Ambrosio L, Fusco N, Pallara D. \emph{Functions of Bounded Variations and Free
  Discontinuity Problems}. Oxford Mathematical Monographs, Oxford Science
  Publications, 2000.

\bibitem{BouFraMar00}
Bourdin B, Francfort GA, Marigo JJ. Numerical experiments in revisited brittle
  fracture. \emph{Journal of the Mechanics and Physics of Solids}  2000;
  \textbf{48}:787--826.

\bibitem{AmbTor92}
Ambrosio L, Tortorelli VM. {On the approximation of free discontinuity
  problems}. \emph{Bollettino dell'Unione Matematica Italiana}  1992;
  \textbf{7}(6-B):105--123.

\bibitem{MumSha89}
Mumford D, Shah J. Optimal approximations by piecewise smooth functions and
  associated variational problems. \emph{Communications on Pure and Applied
  Mathematics}  1989; \textbf{42}:577--685.

\bibitem{Gia05}
Giacomini A. {A}mbrosio-{T}ortorelli approximation of quasi-static evolution of
  brittle fractures. \emph{Calculus of Variations and Partial Differential
  Equations}  2005; \textbf{22}:129--172.

\bibitem{HakKar09}
Hakim V, Karma A. Laws of crack motion and phase-field models of fracture.
  \emph{Journal of the Mechanics and Physics of Solids}  2009;
  \textbf{57}(2):342--368.

\bibitem{LorGod11}
Lorentz E, Godard V. Gradient damage models: toward full-scale computations.
  \emph{Computer Methods in Applied Mechanics and Engineering}  2011;
  \textbf{200}(21--22):1927--1944.

\bibitem{PhaAmoMar11}
Pham K, Amor H, Marigo JJ, Maurini C. Gradient damage models and their use to
  approximate brittle fracture. \emph{International Journal of Damage
  Mechanics}  2011; \textbf{20}(4):618--652.

\bibitem{BouLarRic11}
Bourdin B, Larsen C, Richardson C. A time-discrete model for dynamic fracture
  based on crack regularization. \emph{International Journal of Fracture}
  2011; \textbf{168}(2):133--143.

\bibitem{BorVerSco12}
Borden MJ, Verhoosel CV, Scott MA, Hughes TJ, Landis CM. A phase-field
  description of dynamic brittle fracture. \emph{Computer Methods in Applied
  Mechanics and Engineering}  2012; \textbf{217--220}(0):77--95.

\bibitem{SchWilKuh14}
Schl{\"u}ter A, Willenb{\"u}cher A, Kuhn C, M{\"u}ller R. Phase field
  approximation of dynamic brittle fracture. \emph{Computational Mechanics}
  2014; \textbf{54}(5):1--21.

\bibitem{AbdAri12}
Abdollahi A, Arias I. Phase-field modeling of crack propagation in
  piezoelectric and ferroelectric materials with different electromechanical
  crack conditions. \emph{Journal of the Mechanics and Physics of Solids}
  2012; \textbf{60}(12):2100--2126.

\bibitem{WheWicWol14}
Wheeler MF, Wick T, Wollner W. An augmented-{L}agrangian method for the
  phase-field approach for pressurized fractures. \emph{Computer Methods in
  Applied Mechanics and Engineering}  4 2014; \textbf{271}(0):69--85.

\bibitem{LiPecMil14}
Li B, Peco C, Mill{\'a}n D, Arias I, Arroyo M. Phase-field modeling and
  simulation of fracture in brittle materials with strongly anisotropic surface
  energy. \emph{International Journal for Numerical Methods in Engineering}
  2014; \textbf{102}(3--4):711--727.

\bibitem{DelLanMar07}
Del~Piero G, Lancioni G, March R. A variational model for fracture mechanics:
  numerical experiments. \emph{Journal of the Mechanics and Physics of Solids}
  2007; \textbf{55}:2513--2537.

\bibitem{ClaKna14}
Clayton J, Knap J. A geometrically nonlinear phase field theory of brittle
  fracture. \emph{International Journal of Fracture}  2014;
  \textbf{189}(2):1--10.

\bibitem{HesWei14}
Hesch C, Weinberg K. Thermodynamically consistent algorithms for a
  finite-deformation phase-field approach to fracture. \emph{International
  Journal for Numerical Methods in Engineering}  2014;
  \textbf{99}(12):906--924.

\bibitem{VerBor13}
Verhoosel CV, de~Borst R. A phase-field model for cohesive fracture.
  \emph{International Journal for Numerical Methods in Engineering}  2013;
  \textbf{96}(1):43--62.

\bibitem{AmoMarMau09}
Amor H, Marigo JJ, Maurini C. Regularized formulation of the variational
  brittle fracture with unilateral contact: Numerical experiments.
  \emph{Journal of the Mechanics and Physics of Solids}  2009;
  \textbf{57}(8):1209--1229.

\bibitem{LanRoy09}
Lancioni G, Royer-Carfagni G. The variational approach to fracture mechanics.
  {A} practical application to the {F}rench {P}anth{\'e}on in {P}aris.
  \emph{Journal of Elasticity}  2009; \textbf{95}:1--30.

\bibitem{FreRoy10}
Freddi F, Royer-Carfagni G. Regularized variational theories of fracture: a
  unified approach. \emph{Journal of the Mechanics and Physics of Solids}
  2010; \textbf{58}:1154--1174.

\bibitem{MieWelHof10}
Miehe C, Welschinger F, Hofacker M. Thermodynamically consistent phase-field
  models of fracture: Variational principles and multi-field {FE}
  implementations. \emph{International Journal for Numerical Methods in
  Engineering}  2010; \textbf{83}(10):1273--1311.

\bibitem{AmbGerDe-15}
Ambati M, Gerasimov T, De~Lorenzis L. A review on phase-field models of brittle
  fracture and a new fast hybrid formulation. \emph{Computational Mechanics}
  2015; \textbf{55}(2):383--405.

\bibitem{AmiMilShe14}
Amiri F, Mill{\'a}n D, Shen Y, Rabczuk T, Arroyo M. Phase-field modeling of
  fracture in linear thin shells. \emph{Theoretical and Applied Fracture
  Mechanics}  2014; \textbf{69}:102 -- 109.

\bibitem{AmbLor16}
Ambati M, Lorenzis LD. Phase-field modeling of brittle and ductile fracture in
  shells with isogeometric nurbs-based solid-shell elements. \emph{Computer
  Methods in Applied Mechanics and Engineering}  2016; :--.

\bibitem{LeoBabBou14}
Leon~Baldelli A, Babadjian JF, Bourdin B, Henao D, Maurini C. A variational
  model for fracture and debonding of thin films under in-plane loadings.
  \emph{Journal of the Mechanics and Physics of Solids}  2014;
  \textbf{70}(0):320--348.

\bibitem{PhaMarMau11}
Pham K, Marigo JJ, Maurini C. The issues of the uniqueness and the stability of
  the homogeneous response in uniaxial tests with gradient damage models.
  \emph{Journal of the Mechanics and Physics of Solids}  2011;
  \textbf{59}(6):1163--1190.

\bibitem{PhaMar13b}
Pham K, Marigo JJ. From the onset of damage to rupture: Construction of
  responses with damage localization for a general class of gradient damage
  models. \emph{Continuum Mechanics and Thermodynamics}  2013;
  \textbf{25}(2-4):147--171. Cited By (since 1996)1.

\bibitem{BouMarMau14}
Bourdin B, Marigo JJ, Maurini C, Sicsic P. Morphogenesis and propagation of
  complex cracks induced by thermal shocks. \emph{Physical Review Letters}
  2014; \textbf{112}:014\,301.

\bibitem{Bou07a}
Bourdin B. Numerical implementation of the variational formulation for
  quasi-static brittle fracture. \emph{Interfaces and Free Boundaries}  2007;
  \textbf{9}:411---430.

\bibitem{MieHofWel10a}
Miehe C, Hofacker M, Welschinger F. A phase field model for rate-independent
  crack propagation: Robust algorithmic implementation based on operator
  splits. \emph{Computer Methods in Applied Mechanics and Engineering}  2010;
  \textbf{199}(45{\^a}€``48):2765--2778.

\bibitem{facchinei1998}
Facchinei F, J\'udice J, Soares J. An active set {N}ewton algorithm for
  large-scale nonlinear programs with box constraints. \emph{SIAM Journal on
  Optimization}  1998; \textbf{8}(1):158--186.

\bibitem{ulbrich2011}
Ulbrich M. \emph{Semismooth {N}ewton Methods for Variational Inequalities and
  Constrained Optimization Problems in Function Spaces}, \emph{MOS-SIAM Series
  on Optimization}, vol.~11. SIAM, 2011.

\bibitem{LorBad04}
Lorentz E, Badel P. A new path-following constraint for strain-softening finite
  element simulations. \emph{International Journal for Numerical Methods in
  Engineering}  2004; \textbf{60}(2):499--526.

\bibitem{HeiWheWic15}
Heister T, Wheeler MF, Wick T. A primal-dual active set method and
  predictor-corrector mesh adaptivity for computing fracture propagation using
  a phase-field approach. \emph{Computer Methods in Applied Mechanics and
  Engineering}  2015; \textbf{290}:466--495.

\bibitem{SinVerBor16}
Singh N, Verhoosel C, de~Borst R, van Brummelen E. A fracture-controlled
  path-following technique for phase-field modeling of brittle fracture.
  \emph{Finite Elements in Analysis and Design}  2016; \textbf{113}:14 -- 29.

\bibitem{GerDe-}
Gerasimov T, De~Lorenzis L. A line search assisted monolithic approach for
  phase-field computing of brittle fracture. \emph{Computer Methods in Applied
  Mechanics and Engineering}  in press;
  :--\doi{http://dx.doi.org/10.1016/j.cma.2015.12.017}.
  \urlprefix\url{http://www.sciencedirect.com/science/article/pii/S0045782515004235}.

\bibitem{Neg99}
Negri M. The anisotropy introduced by the mesh in the finite element
  approximation of the {M}umford-{S}hah functional. \emph{Numerical Functional
  Analysis and Optimization}  1999; \textbf{20}:957--982.

\bibitem{BelCos94}
Bellettini G, Coscia A. Discrete approximation of a free discontinuity problem.
  \emph{Numerical Functional Analysis and Optimization}  1994;
  \textbf{15}:105--123.

\bibitem{BurOrtSul10}
Burke S, Ortner C, S\"uli E. An adaptive finite element approximation of a
  variational model of brittle fracture. \emph{SIAM Journal on Numerical
  Analysis}  2010; \textbf{48}(3):980--1012.

\bibitem{ArtForMic15}
Artina M, Fornasier M, Micheletti S, Perotto S. Anisotropic mesh adaptation for
  crack detection in brittle materials. \emph{SIAM Journal on Scientific
  Computing}  2015; \textbf{37}(4):B633--B659.

\bibitem{brune2013}
Brune P, Knepley MG, Smith B, Tu X. Composing scalable nonlinear algebraic
  solvers. \emph{preprint ANL/MCS-P2010-0112}, Argonne National Laboratory
  2013.

\bibitem{ortega1970}
Ortega J, Rheinboldt W. \emph{Iterative Solution of Nonlinear Equations in
  Several Variables}. SIAM, 2000.

\bibitem{young1950}
Young DM. Iterative methods for solving partial difference equations of
  elliptic type. Ph{D} {T}hesis, Harvard University, Massachusetts, USA 1950.

\bibitem{arms1956}
Arms RJ, Gates LD, Zondek B. A method of block iteration. \emph{Journal of the
  Society for Industrial and Applied Mathematics}  1956;
  \textbf{4}(4):220--229.

\bibitem{kahan1958}
Kahan W. Gauss-{S}eidel methods of solving large systems of linear equations.
  Ph{D} {T}hesis, University of Toronto, Ontario, Canada 1958.

\bibitem{ostrowski1954}
Ostrowski AM. On the linear iteration procedures for symmetric matrices.
  \emph{Rendiconti di Matematica e sue Applicazioni}  1954;
  \textbf{5}(14):140--163.

\bibitem{reid1966}
Reid JK. A method for finding the optimum successive over-relaxation parameter.
  \emph{The Computer Journal}  1966; \textbf{9}(2):200--204,
  \doi{10.1093/comjnl/9.2.200}.

\bibitem{munson2000}
Munson TS. Algorithms and environments for complementarity. Ph{D} {T}hesis,
  University of Wisconsin-Madison, Wisconsin, USA 2000.

\bibitem{benson2006}
Benson SJ, Munson TS. Flexible complementarity solvers for large-scale
  applications. \emph{Optimization Methods and Software}  2006;
  \textbf{21}(1):155--168.

\bibitem{fischer1992}
Fischer A. A special {N}ewton-type optimization method. \emph{Optimization}
  1992; \textbf{24}(3--4):269--284.

\bibitem{facchinei1997}
Facchinei F, Soares J. A new merit function for nonlinear complementarity
  problems and a related algorithm. \emph{SIAM Journal on Optimization}  1997;
  \textbf{7}(1):225--247.

\bibitem{hestenes1952}
Hestenes MR, Stiefel E. Methods of {c}onjugate {g}radients for solving linear
  systems. \emph{Journal of Research of the National Bureau of Standards}
  1952; \textbf{49}(6):409--436.

\bibitem{adams2004}
Adams MF, Bayraktar HH, Keaveny TM, Papadopoulos P. Ultrascalable implicit
  finite element analyses in solid mechanics with over a half a billion degrees
  of freedom. \emph{ACM/IEEE Proceedings of SC2004: High Performance Networking
  and Computing}, Pittsburgh, Pennsylvania, 2004.

\bibitem{falgout2006}
Falgout RD. An introduction to algebraic multigrid computing. \emph{Computing
  in Science \& Engineering}  2006; \textbf{8}(6):24--33.

\bibitem{hernandez2005}
Hernandez V, Roman JE, Vidal V. {SLEPc}: A scalable and flexible toolkit for
  the solution of eigenvalue problems. \emph{ACM Transactions on Mathematical
  Software}  2005; \textbf{31}(3):351--362.

\bibitem{vanek1996}
Van\v{e}k P, Mandel J, Brezina M. Algebraic multigrid by smoothed aggregation
  for second and fourth order elliptic problems. \emph{Computing}  1996;
  \textbf{56}(3):179--196.

\bibitem{gee2006}
Gee M, Siefert C, Hu J, Tuminaro R, Sala M. {ML} 5.0 smoothed aggregation
  user's guide. \emph{Technical {R}eport SAND2006-2649}, Sandia National
  Laboratories 2006.

\bibitem{murphy2000}
Murphy MF, Golub GH, Wathen AJ. A note on preconditioning for indefinite linear
  systems. \emph{SIAM Journal on Scientific Computing}  2000;
  \textbf{21}(6):1969--1972.

\bibitem{benzi2005}
Benzi M, Golub GH, Liesen J. Numerical solution of saddle point problems.
  \emph{Acta Numerica}  5 2005; \textbf{14}:1--137.

\bibitem{brown2012}
Brown J, Knepley M, May D, McInnes L, Smith B. Composable linear solvers for
  multiphysics. \emph{2012 11th International Symposium on Parallel and
  Distributed Computing (ISPDC)}, 2012; 55--62.

\bibitem{balay2014}
Balay S, Abhyankar S, Adams MF, Brown J, Brune P, Buschelman K, Eijkhout V,
  Gropp WD, Kaushik D, Knepley MG, \emph{et~al.}. {PETS}c users manual.
  \emph{Technical {R}eport ANL-95/11 - Revision 3.5}, Argonne National
  Laboratory 2014.

\bibitem{paige1975}
Paige C, Saunders M. Solution of sparse indefinite systems of linear equations.
  \emph{SIAM Journal on Numerical Analysis}  1975; \textbf{12}(4):617--629.

\bibitem{HosHsuBou14}
Hossain M, Hsueh CJ, Bourdin B, Bhattacharya K. Effective toughness of
  heterogeneous media. \emph{Journal of the Mechanics and Physics of Solids}
  2014; \textbf{71}(0):15--32.

\bibitem{BahWeiMas88}
Bahr HA, Weiss HJ, Maschke H, Meissner F. Multiple crack propagation in a strip
  caused by thermal shock. \emph{Theoretical and Applied Fracture Mechanics}
  1988; \textbf{10}:219--226.

\bibitem{JiaWuLi12}
Jiang C, Wu X, Li J, Song F, Shao Y, Xu X, Yan P. A study of the mechanism of
  formation and numerical simulations of crack patterns in ceramics subjected
  to thermal shock. \emph{Acta Materialia}  2012; \textbf{60}(11):4540--4550.

\bibitem{SicMarMau14}
Sicsic P, Marigo JJ, Maurini C. Initiation of a periodic array of cracks in the
  thermal shock problem: A gradient damage modeling. \emph{Journal of the
  Mechanics and Physics of Solids}  2014; \textbf{63}(0):256--284.

\bibitem{amestoy2001}
Amestoy PR, Duff IS, Koster J, L'Excellent JY. A fully asynchronous
  multifrontal solver using distributed dynamic scheduling. \emph{SIAM Journal
  on Matrix Analysis and Applications}  2001; \textbf{23}(1):15--41.

\bibitem{eisenstat1996}
Eisenstat S, Walker H. Choosing the forcing terms in an inexact {N}ewton
  method. \emph{SIAM Journal on Scientific Computing}  1996;
  \textbf{17}(1):16--32.

\bibitem{mardal2011}
Mardal KA, Winther R. Preconditioning discretizations of systems of partial
  differential equations. \emph{Numerical Linear Algebra with Applications}
  2011; \textbf{18}(1):1--40.

\end{thebibliography}

\end{document}